\documentclass[reqno,11pt]{amsart}
\usepackage{amsmath}
\usepackage{amsxtra}
\usepackage{amscd}
\usepackage{amsthm}
\usepackage{amsfonts}
\usepackage{amssymb}
\usepackage{eucal}
\usepackage[matrix,arrow,curve]{xy}

\theoremstyle{plain}
\newtheorem{Thm}[subsection]{Theorem}
\newtheorem{Cor}[subsection]{Corollary}
\newtheorem{Lem}[subsection]{Lemma}
\newtheorem{Prop}[subsection]{Proposition}
\newtheorem{Conj}[subsection]{Conjecture}

\theoremstyle{definition}
\newtheorem{Def}[subsection]{Definition}

\theoremstyle{remark}

\newtheorem{Rem}[subsection]{Remark}

\newtheorem{conj}{Conjecture}

\numberwithin{equation}{section}
\renewcommand{\rm}{\normalshape}

\newif\ifShowLabels
\ShowLabelstrue
\newdimen\theight
\def\TeXref#1{%
    \leavevmode\vadjust{\setbox0=\hbox{{\tt
        \quad\quad  {\small \rm #1}}}%
    \theight=\ht0
    \advance\theight by \lineskip
    \kern -\theight \vbox to
    \theight{\rightline{\rlap{\box0}}%
    \vss}%
    }}%

\ShowLabelsfalse

\renewcommand{\sec}[2]{\section{#2}\label{S:#1}%
    \ifShowLabels \TeXref{{S:#1}} \fi}
\newcommand{\ssec}[2]{\subsection{#2}\label{SS:#1}%
    \ifShowLabels \TeXref{{SS:#1}} \fi}

\newcommand{\refs}[1]{Section ~\ref{S:#1}}
\newcommand{\refss}[1]{Section ~\ref{SS:#1}}

\newcommand{\reft}[1]{Theorem ~\ref{T:#1}}
\newcommand{\refl}[1]{Lemma ~\ref{L:#1}}
\newcommand{\refp}[1]{Proposition ~\ref{P:#1}}
\newcommand{\refc}[1]{Corollary ~\ref{C:#1}}

\newcommand{\refe}[1]{\eqref{E:#1}}

\newenvironment{thm}[1]%
    { \begin{Thm} \label{T:#1}  \ifShowLabels \TeXref{T:#1} \fi }%
    { \end{Thm} }

\renewcommand{\th}[1]{\begin{thm}{#1} \sl }
\renewcommand{\eth}{\end{thm} }

\newenvironment{lemma}[1]%
    { \begin{Lem} \label{L:#1}  \ifShowLabels \TeXref{L:#1} \fi }%
    { \end{Lem} }
\newcommand{\lem}[1]{\begin{lemma}{#1} \sl}
\newcommand{\elem}{\end{lemma}}

\newenvironment{propos}[1]%
    { \begin{Prop} \label{P:#1}  \ifShowLabels \TeXref{P:#1} \fi }%
    { \end{Prop} }
\newcommand{\prop}[1]{\begin{propos}{#1}\sl }
\newcommand{\eprop}{\end{propos}}

\newenvironment{corol}[1]%
    { \begin{Cor} \label{C:#1}  \ifShowLabels \TeXref{C:#1} \fi }%
    { \end{Cor} }
\newcommand{\cor}[1]{\begin{corol}{#1} \sl }
\newcommand{\ecor}{\end{corol}}

\newenvironment{defeni}[1]%
    { \begin{Def} \label{D:#1}  \ifShowLabels \TeXref{D:#1} \fi }%
    { \end{Def} }
\newcommand{\defe}[1]{\begin{defeni}{#1} \sl }
\newcommand{\edefe}{\end{defeni}}

\newenvironment{remark}[1]%
    { \begin{Rem} \label{R:#1}  \ifShowLabels \TeXref{R:#1} \fi }%
    { \end{Rem} }
\newcommand{\rem}[1]{\begin{remark}{#1}}
\newcommand{\erem}{\end{remark}}

\newenvironment{conjec}[1]%
    { \begin{Conj} \label{Co:#1}  \ifShowLabels \TeXref{Co:#1} \fi }%
    { \end{Conj} }
\renewcommand{\conj}[1]{\begin{conjec}{#1} \sl }
\newcommand{\econj}{\end{conjec}}

\newcommand{\eq}[1]%
    { \ifShowLabels \TeXref{E:#1} \fi
       \begin{equation} \label{E:#1} }
\newcommand{\eeq}{ \end{equation} }

\newcommand{\prf}{ \begin{proof} }
\newcommand{\epr}{ \end{proof} }

\newcommand{\nc}{\newcommand}

\newcommand{\iso}{\stackrel{\sim}{\longrightarrow}}
\nc{\HC}{{\mathcal{HC}}}
\nc{\on}{\operatorname}
\nc{\BA}{{\mathbb{A}}}
\nc{\BC}{{\mathbb{C}}}
\nc{\BG}{{\mathbb{G}}}
\nc{\BM}{{\mathbb{M}}}
\nc{\BN}{{\mathbb{N}}}
\nc{\BQ}{{\mathbb{Q}}}
\nc{\BP}{{\mathbb{P}}}
\nc{\BR}{{\mathbb{R}}}
\nc{\BZ}{{\mathbb{Z}}}
\nc{\BS}{{\mathbb{S}}}

\nc{\CA}{{\mathcal{A}}}
\nc{\CB}{{\mathcal{B}}}
\nc{\CalC}{{\mathcal C}}
\nc{\CalD}{{\mathcal D}}
\nc{\CE}{{\mathcal{E}}}
\nc{\CF}{{\mathcal{F}}}
\nc{\CG}{{\mathcal{G}}}
\nc{\CH}{{\mathcal{H}}}
\nc{\calJ}{{\mathcal J}}
\nc{\CK}{{\mathcal{K}}}
\nc{\CL}{{\mathcal{L}}}
\nc{\CM}{{\mathcal{M}}}
\nc{\CMM}{{\mathcal{M}^{\operatorname{gen}}_\hbar(-\rho)}}
\nc{\CN}{{\mathcal{N}}}
\nc{\CO}{{\mathcal{O}}}
\nc{\CP}{{\mathcal{P}}}
\nc{\CQ}{{\mathcal{Q}}}
\nc{\CR}{{\mathcal{R}}}
\nc{\CS}{{\mathcal{S}}}
\nc{\CT}{{\mathcal{T}}}
\nc{\CU}{{\mathcal{U}}}
\nc{\CV}{{\mathcal{V}}}
\nc{\CW}{{\mathcal{W}}}
\nc{\CX}{{\mathcal{X}}}
\nc{\CY}{{\mathcal{Y}}}
\nc{\CZ}{{\mathcal{Z}}}

\nc{\gen}{{\operatorname{gen}}}
\nc{\cM}{{\check{\mathcal M}}{}}
\nc{\csM}{{\check{\mathcal A}}{}}
\nc{\obM}{{\overset{\circ}{\mathbf M}}{}}
\nc{\oCA}{{\overset{\circ}{\mathcal A}}{}}
\nc{\obA}{{\overset{\circ}{\mathbf A}}{}}
\nc{\ooM}{{\overset{\circ}{M}}{}}
\nc{\osM}{{\overset{\circ}{\mathsf M}}{}}
\nc{\vM}{{\overset{\bullet}{\mathcal M}}{}}
\nc{\nM}{{\underset{\bullet}{\mathcal M}}{}}
\nc{\obD}{{\overset{\circ}{\mathbf D}}{}}
\nc{\cp}{{\overset{\circ}{\mathbf p}}{}}
\nc{\ofZ}{{\overset{\circ}{\mathfrak Z}}{}}

\nc{\fa}{{\mathfrak{a}}}
\nc{\fb}{{\mathfrak{b}}}
\nc{\fg}{{\mathfrak{g}}}
\nc{\fgl}{{\mathfrak{gl}}}
\nc{\fh}{{\mathfrak{h}}}
\nc{\fj}{{\mathfrak{j}}}
\nc{\fm}{{\mathfrak{m}}}
\nc{\fn}{{\mathfrak{n}}}
\nc{\fu}{{\mathfrak{u}}}
\nc{\fp}{{\mathfrak{p}}}
\nc{\frr}{{\mathfrak{r}}}
\nc{\fs}{{\mathfrak{s}}}
\nc{\ft}{{\mathfrak{t}}}
\nc{\fT}{{\mathfrak{T}}}
\nc{\ofT}{{\overline{\mathfrak T}}}
\nc{\ofS}{{\overline{\mathfrak S}}}
\nc{\fsl}{{\mathfrak{sl}}}
\nc{\hsl}{{\widehat{\mathfrak{sl}}}}
\nc{\hgl}{{\widehat{\mathfrak{gl}}}}
\nc{\hg}{{\widehat{\mathfrak{g}}}}
\nc{\chg}{{\widehat{\mathfrak{g}}}{}^\vee}
\nc{\hn}{{\widehat{\mathfrak{n}}}}
\nc{\chn}{{\widehat{\mathfrak{n}}}{}^\vee}

\nc{\fA}{{\mathfrak{A}}}
\nc{\fB}{{\mathfrak{B}}}
\nc{\fD}{{\mathfrak{D}}}
\nc{\fE}{{\mathfrak{E}}}
\nc{\fF}{{\mathfrak{F}}}
\nc{\fG}{{\mathfrak{G}}}
\nc{\fI}{{\mathfrak{I}}}
\nc{\fJ}{{\mathfrak{J}}}
\nc{\fK}{{\mathfrak{K}}}
\nc{\fL}{{\mathfrak{L}}}
\nc{\fM}{{\mathfrak{M}}}
\nc{\fN}{{\mathfrak{N}}}
\nc{\frP}{{\mathfrak{P}}}
\nc{\fQ}{{\mathfrak Q}}
\nc{\fS}{{\mathfrak S}}
\nc{\fU}{{\mathfrak{U}}}
\nc{\fV}{{\mathfrak{V}}}
\nc{\fZ}{{\mathfrak{Z}}}

\nc{\bb}{{\mathbf{b}}}
\nc{\bc}{{\mathbf{c}}}
\nc{\be}{{\mathbf{e}}}
\nc{\bj}{{\mathbf{j}}}
\nc{\bn}{{\mathbf{n}}}
\nc{\bp}{{\mathbf{p}}}
\nc{\bq}{{\mathbf{q}}}
\nc{\bs}{{\mathbf{s}}}
\nc{\bt}{{\mathbf{t}}}
\nc{\bv}{{\mathbf{v}}}
\nc{\bx}{{\mathbf{x}}}
\nc{\by}{{\mathbf{y}}}
\nc{\bw}{{\mathbf{w}}}
\nc{\bA}{{\mathbf{A}}}
\nc{\bB}{{\mathbf{B}}}
\nc{\bC}{{\mathbf{C}}}
\nc{\bK}{{\mathbf{K}}}
\nc{\bD}{{\mathbf{D}}}
\nc{\bH}{{\mathbf{H}}}
\nc{\bI}{{\mathbf{I}}}
\nc{\bM}{{\mathbf{M}}}
\nc{\bN}{{\mathbf{N}}}
\nc{\bO}{{\mathbf{O}}}
\nc{\bQ}{{\mathbf Q}}
\nc{\bS}{{\mathbf{S}}}
\nc{\bT}{{\mathbf{T}}}
\nc{\bV}{{\mathbf{V}}}
\nc{\bW}{{\mathbf{W}}}
\nc{\bX}{{\mathbf{X}}}
\nc{\bP}{{\mathbf{P}}}
\nc{\bZ}{{\mathbf{Z}}}

\nc{\sA}{{\mathsf{A}}}
\nc{\sB}{{\mathsf{B}}}
\nc{\sC}{{\mathsf{C}}}
\nc{\sD}{{\mathsf{D}}}
\nc{\sF}{{\mathsf{F}}}
\nc{\sK}{{\mathsf{K}}}
\nc{\sM}{{\mathsf{M}}}
\nc{\sO}{{\mathsf{O}}}
\nc{\sQ}{{\mathsf{Q}}}
\nc{\sP}{{\mathsf{P}}}
\nc{\sT}{{\mathsf{T}}}
\nc{\sV}{{\mathsf{V}}}
\nc{\sW}{{\mathsf{W}}}
\nc{\sX}{{\mathsf{X}}}
\nc{\sZ}{{\mathsf{Z}}}
\nc{\sfp}{{\mathsf{p}}}
\nc{\sr}{{\mathsf{r}}}
\nc{\st}{{\mathsf{t}}}
\nc{\sv}{{\mathsf{v}}}
\nc{\sfb}{{\mathsf{b}}}
\nc{\sfc}{{\mathsf{c}}}
\nc{\sd}{{\mathsf{d}}}
\nc{\sg}{{\mathsf{g}}}
\nc{\sfl}{{\mathsf{l}}}

\nc{\BK}{{\bar{K}}}

\nc{\tA}{{\widetilde{\mathbf{A}}}}
\nc{\tB}{{\widetilde{\mathcal{B}}}}
\nc{\tg}{{\widetilde{\mathfrak{g}}}}
\nc{\tG}{{\widetilde{G}}}
\nc{\TM}{{\widetilde{\mathbb{M}}}{}}
\nc{\tO}{{\widetilde{\mathsf{O}}}{}}
\nc{\tU}{{\widetilde{\mathfrak{U}}}{}}
\nc{\TZ}{{\tilde{Z}}}
\nc{\tZ}{\widetilde{Z}{}}
\nc{\tx}{{\tilde{x}}}
\nc{\tbv}{{\tilde{\bv}}}
\nc{\tfP}{{\widetilde{\mathfrak{P}}}{}}
\nc{\tz}{{\tilde{\zeta}}}
\nc{\tmu}{{\tilde{\mu}}}

\nc{\td}{\ddot{\underline{d}}{}}
\nc{\tzeta}{\widetilde{\zeta}{}}
\nc{\hd}{{\widehat{\underline{d}}}}
\nc{\hG}{{\widehat{G}}}
\nc{\hBP}{\widehat{\mathbb P}{}}
\nc{\hQ}{{\widehat{Q}}}
\nc{\hsM}{\widehat{\mathsf M}{}}
\nc{\hfM}{\widehat{\mathfrak M}{}}
\nc{\hCP}{\widehat{\mathcal P}{}}
\nc{\hCR}{\widehat{\mathcal R}{}}
\nc{\hCS}{{\widehat{\mathcal S}}}
\nc{\hfZ}{\widehat{\mathfrak Z}{}}

\nc{\urho}{\underline{\rho}}
\nc{\uB}{\underline{B}}
\nc{\uC}{{\underline{\mathbb{C}}}}
\nc{\ui}{\underline{i}}
\nc{\ofP}{{\overline{\mathfrak{P}}}}

\nc{\hrho}{{\hat{\rho}}}

\nc{\unl}{\underline}
\nc{\ol}{\overline}
\nc{\one}{{\mathbf{1}}}
\nc{\two}{{\mathbf{t}}}

\nc{\Sym}{{\on{Sym}}}
\nc{\Tot}{{\mathop{\operatorname{\rm Tot}}}}
\nc{\Hilb}{{\mathop{\operatorname{\rm Hilb}}}}
\nc{\End}{{\operatorname{End}}}
\nc{\CHom}{{\mathop{\operatorname{{\mathcal{H}}\it om}}}}
\nc{\defi}{{\mathop{\operatorname{\rm def}}}}
\nc{\length}{{\mathop{\operatorname{\rm length}}}}

\nc{\Cliff}{{\mathsf{Cliff}}}
\nc{\Fib}{{\mathsf{Fib}}}
\nc{\Coh}{{\mathsf{Coh}}}
\nc{\FCoh}{{\mathsf{FCoh}}}

\nc{\reg}{{\text{\rm reg}}}
\nc{\oZ}{\overset{\circ}{Z}{}}
\nc{\cplus}{{\mathbf{C}_+}}
\nc{\cminus}{{\mathbf{C}_-}}
\nc{\cthree}{{\mathbf{C}_*}}
\nc{\Qbar}{{\bar{Q}}}
\nc{\Fl}{{{\mathcal F}\ell}}
\nc{\bh}{{\bar{h}}}
\nc{\bOmega}{{\overline{\Omega}}}
\nc\tGr{\widetilde{\Gr}}
\nc{\ul}{\underline}
\nc{\seq}[1]{\stackrel{#1}{\sim}}
\nc\ogu{\overline{G/U}}
\nc\chlam{\check{\lam}}
\newcommand{\Gm}{{\mathbb G}_m}
\nc\St{\operatorname{St}}
\nc{\bLambda}{{\boldsymbol{\Lambda}}}
\nc\uS{\underline{S}}
\nc\QM{\mathcal{QM}}
\nc{\chmu}{\check{\mu}}
\nc{\CHH}{{\CH\!\!\CH}}
\author{Alexander Braverman and Michael Finkelberg}
\title{Twisted zastava and $q$-Whittaker functions}

\begin{document}

\begin{abstract}
We implement the program 
outlined in~\cite[Section~7]{bf11} extending to
the case of non simply laced simple Lie algebras the construction of
solutions of $q$-difference Toda equations from geometry of quasimaps' spaces.
To this end we introduce and study the twisted zastava spaces.
\end{abstract}

\dedicatory{To Valery Lunts on his 60th birthday}
\maketitle 

\sec{Intro}{Introduction}

In this note, we implement the program outlined in~\cite[Section~7]{bf11}
extending to
the case of non simply laced simple Lie algebras the construction of
solutions of $q$-difference Toda equations from geometry of quasimaps' spaces.

\ssec{11}{Semiinfinite Borel-Weil-Bott}
Let $G$ be an almost simple simply connected group over $\BC$ with Lie algebra
$\fg$; we shall denote by $\check\fg$ the Langlands dual algebra of $\fg$.
We fix a Cartan torus and a Borel subgroup $T\subset B\subset G$.
Let also $\CB_{\fg}$ denote its flag variety.
We have $H_2(\CB_{\fg},\BZ)=\Lambda$, the coroot lattice of $\fg$.
We shall denote by $\Lambda_+$ the sub-semigroup
of positive elements in $\Lambda$.

Let $\bC\simeq \BP^1$ denote a (fixed) smooth connected projective curve (over $\BC$) of genus $0$; we are going to fix
a marked point $\infty\in \bC$, and a coordinate $\st$ on $\bC$ such that
$\st(\infty)=0$.
For each $\alpha\in\Lambda_+$ we can consider the space $\CM^{\alpha}_{\fg}$
of maps $\bC\to \CB_{\fg}$ of degree $\alpha$. This is a smooth quasi-projective variety.
It has a compactification $\QM^{\alpha}_{\fg}$ by means of the space of {\em quasi-maps} from $\bC$ to $\CB_{\fg}$ of degree
$\alpha$. Set-theoretically this compactification can be described as follows:
\eq{strat}
\QM^{\alpha}_{\fg}=\bigsqcup\limits_{0\leq \beta\leq \alpha}\CM^{\beta}_{\fg}\times\Sym^{\alpha-\beta}(\bC)
\end{equation}
where $\Sym^{\alpha-\beta}(\bC)$ stands for the space of ``colored divisors" of the form
$\sum \gamma_i x_i$ where $x_i\in \bC$, $\gamma_i\in \Lambda_+$ and $\sum \gamma_i=\alpha-\beta$. In particular, for $\beta\geq\alpha$ we have an embedding
$\varphi_{\alpha,\beta}:\ \QM^{\alpha}_{\fg}\hookrightarrow\QM^\beta_{\fg}$ adding
defect at the point $0\in\bC$ (such that $\st(0)=\infty$). The union of all
$\QM^{\alpha}_{\fg}$ is an ind-projective scheme $\fQ_\fg$. To each weight
$\check\lambda\in X^*(T)$ of $G$ one associates a line bundle
$\CO(\check\lambda)$ on $\fQ_\fg$.

Recall the notion of (global) Weyl modules $\CW(\check\lambda)$ over the
current algebra $\fg[\st]$ (see e.g.~\cite{cfk}). The following version of
the Borel-Weil-Bott theorem was proved in~\cite{bf12} in case $\fg$ is
simply-laced. First, the higher cohomology $H^{>0}(\fQ_\fg,\CO(\check\lambda))$
vanish identically. Second, in case $\check\lambda$ is {\em not} a dominant
weight, the global sections $H^0(\fQ_\fg,\CO(\check\lambda))$ vanish as well.
Third, in case $\check\lambda$ is a dominant weight, the global sections
$H^0(\fQ_\fg,\CO(\check\lambda))$ are isomorphic to the {\em dual} global
Weyl module $\CW(\check\lambda)^\vee$. In the last~\refs{nono} of the present
note we extend the
Borel-Weil-Bott theorem to the case of arbitrary simple $\fg$, and also prove
that the schemes $\QM^{\alpha}_{\fg}$ have rational singularities.

\ssec{12}{The $q$-Whittaker functions}
Let $\check{G}$ denote the Langlands dual
group of $G$ with its maximal torus $\check{T}$. Let $W$ be the Weyl group
of $(G,T)$. We recall the notion of
$q$-Whittaker functions $\Psi_{\check\lambda}(q,z)\colon W$-invariant polynomials
in $z\in T$ with coefficients in the field $\BC(q)$ of rational functions in 
a variable $q$ ($\check\lambda\in X^*(T)^+$ a dominant weight of $G$).
The definition of $\Psi_{\check\lambda}(q,z)$ is as follows.
In \cite{Et} and \cite{Sev} the authors define (by adapting the so called
Kostant-Whittaker reduction to the case of quantum groups) a homomorphism
$\CM:\ \BC[T]^W\to \End_{\BC(q)}\BC(q)[\check T]$ called the quantum difference
Toda integrable system associated with $\check G$.
For each $f\in \BC[T]^W$ the operator $\CM_f:=\CM(f)$ is indeed a difference
operator: it is a $\BC(q)$-linear combination of
shift operators $\bT_{\check\beta}$ where $\check\beta\in X^*(T)$ and
$$
\bT_{\check\beta}(F(q,x))=F(q,q^{\check\beta}x),\ x\in\check{T}.
$$
In particular, the above operators can be restricted to operators acting in
the space of functions on the lattice $X^*(T)$ by means
of the embedding $X^*(T)\hookrightarrow \check T$ sending every
$\check\lambda$ to $q^{\check\lambda}$. More precisely, we have a restriction
morphism $\on{res}\colon \BC(q)[\check T]\to\BC(q)[X^*(T)],\ \on{res}F(q,\check\lambda):=F(q,q^{\check\lambda})$, and a unique operator $\bT_{\check\beta}\colon
\BC(q)[X^*(T)]\to\BC(q)[X^*(T)]$ such that 
$\bT_{\check\beta}\on{res}=\on{res}\bT_{\check\beta}$. Namely,
$\bT_{\check\beta}(F(q,\check\lambda))=F(q,\check\beta+\check\lambda)$.
For any $f\in \BC[T]^W$ we shall denote the corresponding operator on
$\BC(q)[X^*(T)]$ by $\CM_f^{\on{lat}}$.

There exists a collection of $\BC(q)$-valued
polynomials\footnote{Conjecturally, such a collection is unique. We do not use
the uniqueness property in the present note. The existence follows e.g.\ from 
the results of the present note.} 
of  $\Psi_{\check\lambda}(q,z),\ \check\lambda\in X^*(T)$, on $T$
satisfying the following properties:

a) $\Psi_{\check\lambda}(q,z)=0$ if $\check\lambda$ is not dominant.

b) $\Psi_{0}(q,z)=1$.

c) Let us consider all the functions $\Psi_{\check\lambda}(q,z)$ as one function
$\Psi(q,z):\ X^*(T)\to \BC(q)$ depending
on $z\in  T$. Then for every $f\in \BC[T]^W$ we have
$$
\CM^{\on{lat}}_f(\Psi(q,z))=f(z)\Psi(q,z).
$$

There exists another definition of the $q$-Toda system using double affine
Hecke algebras, studied for example in~\cite{c}. To be more specific, we
restrict ourselves here to  the double affine Hecke algebras of symmetric
type in terminology of~\cite{h}. Since it is not clear to us how to prove
{\em apriori} that the definition of $q$-Toda from~\cite{c} coincides with
the definitions from~\cite{Et} and~\cite{Sev}, we shall denote the
$q$-difference operators from~\cite{c} by $\CM'_f$. Similarly we shall
denote by $(\CM_f^{\on{lat}})'$ their ``lattice'' version. We shall denote the
corresponding polynomials by $\Psi'_{\check\lambda}(q,z)$, so that
$(\CM^{\on{lat}}_f)'(\Psi'(q,z))=f(z)\Psi'(q,z)$. Their existence follows e.g.\
from the results of~\cite{i}.

\ssec{13}{Characters of twisted Weyl modules}
In case $\fg$ is simply laced, it was proved in~\cite{bf12} that
$\Psi_{\check\lambda}(q,z)$ coincides with the character of the global Weyl module
$\CW(\check\lambda)$ over $\fg[\st]\rtimes\BC^*$; and it was explained
in~Section~1.4 of {\em loc. cit.}
that such an equality does not hold in case of non simply laced $\fg$.
In the non simply laced case we use the following remedy.
We realize $\check\fg$ as a {\em folding} of a simple
simply laced Lie algebra
$\check\fg'$, i.e. as invariants of an outer automorphism $\sigma$ of
$\check\fg'$ preserving a Cartan subalgebra $\check\ft'\subset\check\fg'$
and acting on the root system of $(\check\fg',\check\ft')$.
In particular,
$\sigma$ gives rise to the same named automorphism of the Langlands dual
Lie algebras $\fg'\supset\ft'$ (note that say, in case $\fg$ is of type
$B_n$, $\fg'$ is of type $A_{2n-1}$, while for $\fg$ of type
$C_n$, $\fg'$ is of type $D_{n+1}$; in particular, $\fg\not\subset\fg'$).
Let $d$ stand for the order of $\sigma$.
We choose a primitive root of unity $\zeta$ of order $d$.
We consider an automorphism $\varsigma$ of
$\fg'[\st]$ defined as the composition of two automorphisms:
a) $\sigma$ of $\fg'$;
b) $\st\mapsto\zeta\st$ of $\BC[\st]$. The subalgebra of invariants
$\fg'[\st]^\varsigma$ is the twisted current algebra. The corresponding global
twisted Weyl modules $\CW^{\on{twisted}}(\check\lambda)$ over
$\fg'[\st]^\varsigma\rtimes\BC^*$
(still numbered by the dominant $\fg$-weights $\check\lambda\in X^*(T)^+$)
were introduced in~\cite{cik}.

In~\refs{fer} of the present note we prove that the $q$-Whittaker function
$\Psi_{\check\lambda}(q,z)$ coincides with the character of the global
twisted Weyl module
$\CW^{\on{twisted}}(\check\lambda)$ over $\fg'[\st]^\varsigma\rtimes\BC^*$.
The relation between the global and local twisted Weyl modules established
in~\cite{cik} then implies the following positivity property of
$\Psi_{\check\lambda}(q,z)$. Let $d_i=1$ (resp. $d_i=d$) for a short
(resp. long) simple coroot $\alpha_i$ of $\fg$. For $i\in I$: the set of
simple coroots of $\fg$, we set $q_i:=q^{d_i}$. We set
$\hat{\Psi}_{{\check\lambda}}(q,z):=\Psi_{{\check\lambda}}(q,z)\cdot{\prod\limits_{i\in I}
\prod\limits_{r=1}^{\langle\alpha_i,\check\lambda\rangle}(1-q_i^r)}$. Then
$\hat{\Psi}_{{\check\lambda}}(q,z)$ is a polynomial in $z,q$ with nonnegative
integral coefficients. Namely, $\hat{\Psi}_{{\check\lambda}}(q,z)$ is the
character of the local twisted Weyl module.

In fact, the above results are known if one replaces $\hat{\Psi}_{{\check\lambda}}(q,z)$ with the
polynomials $\hat{\Psi}'_{{\check\lambda}}(q,z):=\Psi'_{{\check\lambda}}(q,z)\cdot{\prod\limits_{i\in I}
\prod\limits_{r=1}^{\langle\alpha_i,\check\lambda\rangle}(1-q_i^r)}$ (these are often called $q$-Hermite polynomials
in the literature). Namely, the above local twisted Weyl modules
 coincide by~\cite{fok} with the level
one Demazure module $D^{\on{twisted}}(\check\lambda)$ over
$\fg'[\st]^\varsigma\rtimes\BC^*$. Now the characters of level one Demazure
modules over dual untwisted affine Lie algebras were proved in~\cite{i} to
coincide with the $q$-Hermite polynomials
$\hat{\Psi}'_{{\check\lambda}}(q,z)$.
Thus we obtain the following corollary:
\cor{psi=psi}
We have ${\Psi}_{{\check\lambda}}(q,z)={\Psi}'_{{\check\lambda}}(q,z)$.
Hence the $q$-Toda systems of~\cite{Et},~\cite{Sev} and of~\cite{c} are
equivalent (they generate the same commutative subalgebras in the $q$-difference
operators' rings).
\ecor
Let us note that the above proof of \refc{psi=psi} is very roundabout.
 It would be nice to find
a more direct argument.
\ssec{14}{Twisted quasimaps}
Our proof of the properties~\refss{12}(a,b,c) of the characters of the global
twisted
Weyl modules uses a twisted version of the semiinfinite Borel-Weil-Bott theorem
of~\refss{11}. Namely, the automorphism $\varsigma$ of $\fg'[\st]$ gives rise
to the same named automorphism $\varsigma$ of the ind-projective scheme
$\fQ_{\fg'}$ of~\refss{11}. Its fixed point subscheme is denoted by $\fQ$.
To each weight
$\check\lambda\in X^*(T)$ of $G$ one associates a line bundle
$\CO(\check\lambda)$ on $\fQ$. As in~\refss{11}, we have
$H^{>0}(\fQ,\CO(\check\lambda))=0$, while
$H^0(\fQ,\CO(\check\lambda))=\CW^{\on{twisted}}(\check\lambda)^\vee$.

Now the $q$-difference equations of~\refss{12}c) for the characters of
$H^0(\fQ,\CO(\check\lambda))$ are proved following the strategy
of~\cite{bf11},~\cite{bf12} provided we know some favourable geometric
properties of the finite-type pieces $\QM^\alpha\subset\fQ$ (twisted quasimaps'
spaces: the fixed point sets of the automorphism $\varsigma$ of certain
quasimaps' spaces $\QM^\beta_{\fg'}$) and their local (based) analogues:
twisted zastava spaces $Z^\alpha$. The verification of these properties
occupies the bulk of the present note, namely the central~\refs{twi qua}.
Some properties, like irreducibility and normality of $Z^\alpha$ are proved
similarly to their classical (nontwisted) counterparts, by reduction to the
known properties of the twisted affine Grassmannian of $\fg'$.
Some other, like the Cartier property of the (reduced) boundary, turn out
harder to prove. Very roughly speaking, our derivation of the $q$-difference
equations for the characters of twisted global Weyl modules from the 
semiinfinite Borel-Weil-Bott theorem is parallel to the derivation of the 
Weyl character formula from the classical BWB theorem via localization to the 
torus-fixed points.
Finally, note that the previous results 
of~\cite{bf11},~\cite{bf12} are formally contained in the results of the 
present note in case of trivial folding when $\sigma=1$ and $d=1$.

\ssec{ack}{Acknowledgments} We are grateful to A.~Kuznetsov,
S.~Kov\'acs, Yu.~Prokhorov, L.~Rybnikov, J.~Shiraishi, and X.~Zhu for their
patient explanations.
M.F.\ was partially supported by a subsidy granted to the HSE by
the Government of the Russian Federation for the implementation of the Global
Competitiveness Program.

\sec{globa}{Setup and notations}

\ssec{fold}{Root systems and foldings}
Let $\check\fg$ be a simple Lie algebra
with the corresponding adjoint Lie group $\check G$. Let $\check T$ be
a Cartan torus of $\check G$. We choose a Borel subgroup $\check{B}\supset
\check T$. It defines the set of simple roots $\{\alpha_i,\ i\in I\}$.
Let $G\supset T$ be the Langlands dual groups. We define an isomorphism
$\alpha\mapsto\alpha^*$ from the root lattice of $(\check{G},\check{T})$ to
the root lattice of $(G,T)$ in the basis of simple roots as follows:
$\alpha_i^*:=\check{\alpha}_i$ (the corresponding simple coroot).
For two elements $\alpha,\beta$ of the root lattice of $(\check{G},\check{T})$
we say $\beta\leq\alpha$
if $\alpha-\beta$ is a nonnegative linear combination of
$\{\alpha_i,\ i\in I\}$.
For such $\alpha$ we denote by $z^{\alpha^*}$ the
corresponding character of $T$. As usually, $q$ stands for the
identity character of $\BG_m$. We set $d_i=\frac{(\alpha_i,\alpha_i)}{2}$,
and $q_i=q^{d_i}$.

We realize $\check\fg$ as a {\em folding} of a simple
simply laced Lie algebra
$\check\fg'$, i.e. as invariants of an outer automorphism $\sigma$ of
$\check\fg'$ preserving a Cartan subalgebra $\check\ft'\subset\check\fg'$
and a Borel subalgebra $\check\fb'\supset\check\ft'$,
and acting on the root system of $(\check\fg',\check\ft')$.
Note that the unfolding $(\check\fg'\supset\check\fb'\supset\check\ft',\sigma)$
is defined uniquely up to an isomorphism.
In particular,
$\sigma$ gives rise to the same named automorphism of the Langlands dual
Lie algebras $\fg'\supset\ft'$. We choose a $\sigma$-invariant Borel
subalgebra $\ft'\subset\fb'\subset\fg'$ such that $\fb=(\fb')^\sigma$.
The corresponding set of simple roots is denoted by $I'$.
We denote by $\Xi$ the finite cyclic group generated by $\sigma$.
We set $d:=|\Xi|$. Note that $d_i\in\{1,d\}$.
Let $G'\supset T'$ denote the simply
connected Lie group and its Cartan torus with Lie algebras $\fg'\supset\ft'$.
The {\em coinvariants} $X_*(T')_\sigma$ of $\sigma$ on the coroot
lattice $X_*(T')$ of $(\fg',\ft')$ coincide with the root lattice
of $\check\fg$ and with the coroot lattice of $\fg$. 
We have an injective map $a:\ X_*(T')_\sigma\to X_*(T')^\sigma$
from coinvariants to invariants defined as follows: given a coinvariant
$\alpha$ with a representative $\tilde\alpha\in X_*(T')$ we set
$a(\alpha):=\sum_{\xi\in\Xi}\xi(\tilde\alpha)$.
The Weyl group $W$ of $G\supset T$ coincides with the invariants
$(W')^\sigma$ of $\sigma$ on the Weyl group $W'$ of $G'\supset T'$.
We fix a primitive root of unity $\zeta$ of order $d$.
We set $\CK=\BC((\st))\supset\CO=\BC[[\st]]$. We set $\bt:=\st^{-1}$.



\ssec{QQ}{Ind-scheme $\fQ$}
We denote by $\on{Gr}$ the twisted affine Grassmannian
$G'(\CK)^\varsigma/G'(\CO)^\varsigma$: an ind-proper ind-scheme of ind-finite
type, see~\cite{pr},~\cite{z}.
We consider the projective
line $\bC$ with coordinate $\bt$, and with points $0=0_\bC,\ \infty=\infty_\bC$
such that $\bt(0_\bC)=0,\bt(\infty_\bC)=\infty$. We recall the setup
of~\cite[Section~2]{bf12} with $\fg'$ (resp. $\bt$) playing the role of $\fg$
(resp. $t$) of {\em loc. cit.} In particular, $R=\BC[[t^{-1}]]$
(resp. $F=\BC((t^{-1}))$) of {\em loc. cit.} is our $\CO=\BC[[\st]]$
(resp. $\CK=\BC((\st))$). Furthermore, $\Lambda_+$ of {\em loc. cit.} is the
cone in $X_*(T')$ generated over $\BN$ by the simple coroots, while
$\Lambda^\vee_+$ of {\em loc. cit.} is the cone in $X^*(T')$ generated over
$\BN$ by the fundamental weights. Given $\gamma\in\Lambda_+$, we consider
the quasimaps' space $\QM^\gamma_{\fg'}$.

Recall the notations of~\refss{fold}.
We consider the cone $Y_+\subset Y=X_*(T')_\sigma$ generated over $\BN$ by
the classes of simple coroots of $\fg$.
Given $\alpha\in Y_+$, we consider an automorphism $\varsigma$ of
$\QM^{a(\alpha)}_{\fg'}$ defined as the composition of two automorphisms:
a) $\sigma$ (arising from the same named automorphism of $G'$);
b) $\bt\mapsto\zeta^{-1}\bt$.
We define $\QM^\alpha$ as the fixed point set $(\QM^{a(\alpha)}_{\fg'})^\varsigma$
equipped with the structure of reduced closed subscheme of $\QM^{a(\alpha)}_{\fg'}$.

For $\beta\geq\alpha\in Y_+$ (that is, $\beta-\alpha\in Y_+$), we consider
the closed embedding $\varphi_{\alpha,\beta}:\ \QM^\alpha\hookrightarrow\QM^\beta$
adding the defect $a(\beta-\alpha)\cdot0$ at the point $0\in\bC$.
The direct limit of this system is denoted by $\fQ$.

\ssec{QQQ}{Infinite type scheme $\bQ$}
We fix a collection of highest weight vectors
$\sv_{\check\lambda}\in V_{\check\lambda},\ \check\lambda\in\Lambda^\vee_+\subset
X^*(T')$, satisfying
the Pl\"ucker equations. We denote by $\sigma:\ V_{\check\lambda}\to
V_{\sigma(\check\lambda)}$ a unique isomorphism taking $\sv_{\check\lambda}$
to $\sv_{\sigma(\check\lambda)}$ and intertwining $\sigma:\ G'\to G'$.
We denote by $\widehat\bQ$ the {\em reduced} infinite type scheme whose 
$\BC$-points are the collections of {\em nonzero} vectors
$v_{\check\lambda}(\bt)\in V_{\check\lambda}\otimes\BC[[\bt^{-1}]],\
\check\lambda\in\Lambda^\vee_+$, satisfying the Pl\"ucker relations and
the equation
$\sigma(v_{\check\lambda})(\zeta^{-1}\bt)=v_{\sigma(\check\lambda)}(\bt)$.
It is equipped with a free action of $T=(T')^\sigma$: if we view an element
of $T$ as a $\sigma$-invariant element $h\in(T')^\sigma$, then
$h(v_{\check\lambda}(\bt))=\check\lambda(h)v_{\check\lambda}(\bt)$.
The quotient scheme $\bQ=\widehat\bQ/T$ is a closed subscheme in
$\prod_{i\in I'}\BP(V_{\check\omega_i}\otimes\BC[[\bt^{-1}]])$ where $\check\omega_i$
is a fundamental weight of $\fg'$. Any weight
$\check\lambda\in\Lambda^\vee_\sigma=
X^*(T')_\sigma=X$ gives rise to a line bundle $\CO_{\check\lambda}$ on $\bQ$.

The construction of~\cite[2.3]{bf12} gives rise to the closed embedding
$\fQ\hookrightarrow\bQ$.


Finally, recall that the restriction of characters gives rise to a canonical
isomorphism $X=X^*(T')_\sigma\iso X^*(T)$. The $T$-torsor
$\widehat\bQ\to\bQ$ defines, for any $\check\lambda\in X$,
a line bundle $\CO(\check\lambda)$ on $\bQ$. Same notation for its
restriction to $\fQ$.

\ssec{twist z}{Twisted zastava}
The twisted quasimaps' space $\QM^\alpha=(\QM^{a(\alpha)}_{\fg'})^\varsigma$
has an open dense subvariety $'\QM^\alpha$ formed
by the quasimaps without defect at $\infty\in\bC$. We have an evaluation
morphism $ev_\infty:\ '\QM^\alpha\to\CB:=\CB_{\fg'}^\sigma=(G'/B')^\sigma$.
We define the twisted zastava space
$Z^\alpha:=ev_\infty^{-1}(\fb_-)=(Z^{a(\alpha)}_{\fg'})^\varsigma$.
Recall the factorization morphism
$\pi:\ Z^{a(\alpha)}_{\fg'}\to\BA^{a(\alpha)}:=(\bC-\infty)^{a(\alpha)}$.
We consider an automorphism $\varsigma$ of the coloured divisors' space
$\BA^{a(\alpha)}$ defined as the
composition of two automorphisms: a) $\sigma$ on the set of colours;
b) $\bt\mapsto\zeta^{-1}\bt$ on $\BA^1$. We have
$(\BA^{a(\alpha)})^\varsigma=\BA^\alpha$; a few words about the meaning of the
notation $\BA^\alpha$ are in order. Let $\alpha=\sum_{i\in I}a_i\alpha_i$
where $I=I'/\Xi$ (the orbits of the cyclic group generated by
$\sigma)=I_0\sqcup I_1$ where $I_0$ consists of one-point-orbits (fixed
points), while $I_1$ consists of free orbits (so that $\alpha_i$ is a
long (resp. short) simple root of $(\check G,\check T)$ if $i\in I_0$
(resp. $i\in I_1$)). Then $\BA^\alpha=\prod_{i\in I_1}(\bC-\infty)^{(a_i)}\times
\prod_{i\in I_0}((\bC-\infty)/(\bt\mapsto\zeta^{-1}\bt))^{(a_i)}$.
Note that $(\bC-\infty)/(\bt\mapsto\zeta^{-1}\bt)\simeq\BA^1$ with coordinate
$\bt^d$ (where $d=|\Xi|$, see~\refss{fold}). In particular,
the diagonal stratification of $\BA^{a(\alpha)}$ induces a {\em quasidiagonal}
stratification of $\BA^\alpha$: a point $\unl{w}\in\BA^\alpha$ lies on a
quasidiagonal if either of the following holds: a) $w_{i,r}=w_{j,s}$ for
$i,j\in I_0$ or $i,j\in I_1$ (and $1\leq r\leq a_i,\ 1\leq s\leq a_j$);
b) $w_{i,r}=w_{j,s}^{d}$ for $i\in I_0,\ j\in I_1$.

Now $\pi$ commutes with $\varsigma$, so that the following diagram commutes:
\eq{pisigma}
\begin{CD}
Z^\alpha @>>> Z^{a(\alpha)}_{\fg'}\\
@VVV         @V\pi VV\\
\BA^\alpha @>>> \BA^{a(\alpha)}
\end{CD}
\end{equation}
We will denote the left vertical arrow by $\pi$ as well. The commutativity
of the diagram~\refe{pisigma} implies that the factorization property holds
for $\pi:\ Z^\alpha\to\BA^\alpha$.

\ssec{exa}{An example}
We take $\fg'=\mathfrak{sl}(4)\supset(\fg')^\sigma=\mathfrak{sp}(4)$ 
(the invariants of the outer automorphism). We denote the simple coroots of 
$\fg\simeq\mathfrak{sp}(4)$ by
$\alpha_1,\alpha_2$, and the simple coroots of $\fg'$ by
$\beta_1,\beta_2,\beta_3$, so that $a(\alpha_1)=\beta_1+\beta_3$, and
$a(\alpha_2)=2\beta_2$. We will exhibit an explicit system of equations
defining the twisted zastava $Z^\alpha$ for $\alpha=\alpha_1+\alpha_2$.

To this end recall the fundamental representations of $\fg':\ V=V_{\check\omega_1}$
with a base $v_1,v_2,v_3,v_4;\ \Lambda^2V=V_{\check\omega_2}$ with a base
$v_{ij}:=v_i\wedge v_j,\ 1\leq i<j\leq4$, and finally $\Lambda^3V=V_{\check\omega_3}$
with a base $v_{ijk}:=v_i\wedge v_j\wedge v_k,\ 1\leq i<j<k\leq4$.
The involutive outer automorphism $\sigma$ takes $V$ to $\Lambda^3V$, and
$\Lambda^2V$ to itself; its action in the above bases is as follows:
$v_1\mapsto v_{123},\ v_2\mapsto v_{124},\ v_3\mapsto v_{134},\ v_4\mapsto v_{234};\
v_{12}\mapsto v_{12},\ v_{13}\mapsto v_{13},\ v_{24}\mapsto v_{24},\
v_{34}\mapsto v_{34},\ v_{14}\mapsto -v_{23},\ v_{23}\mapsto -v_{14}$.

Zastava space $Z^{(1,2,1)}_{\mathfrak{sl}(4)}$ is formed by the collections of
$V_{\check\omega_i}$-valued polynomials of the form
$(\bt-a_1)v_1+a_2v_2+a_3v_3+a_4v_4,\
(\bt-a_{123})v_{123}+a_{124}v_{124}+a_{134}v_{134}+a_{234}v_{234},\
(\bt^2+b_{12}\bt-a_{12})v_{12}+(b_{13}\bt+a_{13})v_{13}+(b_{24}\bt+a_{24})v_{24}+
(b_{34}\bt+a_{34})v_{34}+(b_{14}\bt+a_{14})v_{14}+(b_{23}\bt+a_{23})v_{23}$
subject to the Pl\"ucker relations to be specified below.
The twisted zastava space $Z^{(1,1)}\subset Z^{(1,2,1)}_{\mathfrak{sl}(4)}$ is cut
out by the following invariance conditions:
$a_{123}=-a_1,\ a_{124}=-a_2,\ a_{134}=-a_3,\ a_{234}=-a_4,\
b_{12}=b_{13}=b_{24}=b_{34}=0,\ b_{23}=b_{14},\ a_{23}=-a_{14}$.

When writing down the Pl\"ucker relations explicitly we will make use of
the above invariance conditions to simplify the resulting equations.
First, the $\mathfrak{sl}(4)$-invariant projection $V\otimes\Lambda^3V\to\BC$
must annihilate our polynomials, that is $a_{234}-a_4=0$ and
$a_3a_{124}+a_4a_{123}-a_1a_{234}-a_2a_{134}=0$. Substituting the invariance
conditions we get $a_4=a_{234}=0$.
Second, the $\mathfrak{sl}(4)$-invariant projection
$\Lambda^2V\otimes\Lambda^2V\to\BC$ must annihilate our polynomials, that is
$a_{34}+b_{14}b_{23}=0,\ b_{14}a_{23}+b_{23}a_{14}=0,\
a_{14}a_{23}-a_{12}a_{34}-a_{13}a_{24}=0$.
Third, the $\mathfrak{sl}(4)$-invariant projection
$V\otimes\Lambda^2V\to\Lambda^3V$ must annihilate our polynomials, that is
$a_3+b_{23}=0,\ a_4=0;\ a_{24}-a_2b_{14}=0,\ a_{34}-a_3b_{14}=0,\ a_4b_{23}=0,\
a_{23}-a_1b_{23}=0;\ a_1a_{23}+a_2a_{13}+a_3a_{12}=0,\ a_1a_{24}+a_2a_{14}+a_4a_{12}=0,\
a_1a_{34}+a_3a_{14}-a_4a_{13}=0,\ a_2a_{34}-a_3a_{24}+a_4a_{23}=0$.

All in all, we have $a_4=0,\ b_{23}=b_{14}=-a_3,\ a_{23}=-a_{14}$; substituting
for $a_{34},a_{24},a_{14}$ their values from the third group of equations,
we are left with the variables $a_1,a_2,a_3,a_{12},a_{13}$ satisfying the
{\em single} equation $a_3(a_1^2-a_{12})=a_2a_{13}$.
The factorization projection $\pi:\ Z^\alpha\to\BA^\alpha$ sends
$(a_1,a_2,a_3,a_{12},a_{13})$ to $(a_1,a_{12})$. The boundary $\partial Z^\alpha=
Z^\alpha\setminus\overset{\circ}{Z}{}^\alpha$ is given by a single equation
$a_3=0$.

\ssec{anexa}{Another example}
We take $\fg'=\mathfrak{so}(8)\supset(\fg')^\sigma$ 
(the invariants of the outer automorphism of order 3). 
We denote the simple coroots of $\fg$ of type $G_2$ by
$\alpha_1,\alpha_2$, and the simple coroots of $\fg'$ by
$\beta_1,\beta_2,\beta_3,\beta_4$, so that $a(\alpha_1)=3\beta_1$, and
$a(\alpha_2)=\beta_2+\beta_3+\beta_4$. 
We will exhibit an explicit system of equations
defining the twisted zastava $Z^\alpha$ for $\alpha=\alpha_1+\alpha_2$,
so that $\beta=a(\alpha)=3\beta_1+\beta_2+\beta_3+\beta_4$.

Zastava space $Z^\beta_{\mathfrak{so}(8)}$ is formed by the collections of
$V_{\check\omega_i}$-valued polynomials of the form
$(\bt^3-f''\bt^2-f'\bt-f)v_{\check\omega_1}+
(e''\bt^2+e'\bt+e)v_{\check\omega_1-\check\alpha_1}+\ldots,\
(\bt-a)v_{\check\omega_2}+bv_{\check\omega_2-\check\alpha_2}+
cv_{\check\omega_2-\check\alpha_2-\check\alpha_1}+\ldots$, and so on.
One can show that the invariance conditions together with Pl\"ucker equations
boil down to a single equation cutting out $Z^\alpha$ in $\BA^5$ with
coordinates $(a,b,c,e,f)$, namely
$c(a^3-f)=be$. The factorization projection $\pi\colon Z^\alpha\to\BA^\alpha$
sends $(a,b,c,e,f)$ to $(a,f)$. The boundary $\partial Z^\alpha$ is given by
a single equation $c=0$.

In effect, let $Y^\alpha\stackrel{\varpi}{\longrightarrow}\BA^\alpha$ denote
the above hypersurface $c(a^3-f)=be$ in $\BA^5$, and its projection
$(a,b,c,e,f)\mapsto(a,f)$ to $\BA^2$. Then the open subvarieties
$\pi^{-1}(\BA^2\setminus\{(0,0)\})\subset Z^\alpha$ and 
$\varpi^{-1}(\BA^2\setminus\{(0,0)\})\subset Y^\alpha$ are isomorphic 
by~\refl{F1} below and e.g.~\cite[5.6]{bdf}. This isomorphism extends to
$Y^\alpha\simeq Z^\alpha$ due to normality of $Z^\alpha$~(\refp{norm} below).

\sec{twi qua}{Geometric properties of twisted quasimaps}

\ssec{diafi}{Quasidiagonal fibers}
The factorization property of $\pi:\ Z^\alpha\to\BA^\alpha$ implies that
in order to describe the fibers of
$\pi$ it suffices to describe the quasidiagonal fibers
$\CF_0^\alpha:=\pi^{-1}(\alpha\cdot0)$, and $\CF_1^\alpha:=\pi^{-1}(\alpha\cdot1)$
(note that $\CF_1^\alpha$ is isomorphic to 
$\pi^{-1}(\alpha_0\cdot c^{d}+\alpha_1\cdot c)$ for any
$c\ne0$, where $\alpha_0:=\sum_{i\in I_0}a_i\alpha_i$, and
$\alpha_1:=\sum_{i\in I_1}a_i\alpha_i$).
Recall that the diagonal fiber $\pi^{-1}(\gamma\cdot c)\subset Z^\gamma_{\fg'}$
is denoted by $\CF^\gamma_{\fg'}$
(these fibers are all canonically isomorphic for various choices of
$c\in\BA^1$); it is equidimensional of dimension $|\gamma|$. Let us choose
a decomposition $a(\alpha)=\sum_{\xi\in\Xi}\xi(\tilde\alpha)$ as in~\refss{fold}
for $\tilde\alpha\in\Lambda_+\subset X_*(T')$.

\lem{F1}
a) $\CF_1^\alpha\supset\CF^{\tilde\alpha}_{\fg'}$;

b) $\CF_1^\alpha=\left(\bigcup_{\tilde\alpha}\CF^{\tilde\alpha}_{\fg'}\right)/\Xi$ 
(the union over
all the choices of $\tilde\alpha\in\Lambda_+\subset X_*(T')$ such that
$a(\alpha)=\sum_{\xi\in\Xi}\xi(\tilde\alpha)$);

c) In particular, $\dim\CF_1^\alpha=|\alpha|$.
\elem

\prf
Clear.
\epr

In order to describe the (quasi)diagonal fiber $\CF^\alpha_0$ we need the twisted
affine Grassmannian $\on{Gr}=G'(\CK)^\varsigma/G'(\CO)^\varsigma$ of~\refss{QQ}.
The $T$-fixed points of $\on{Gr}$ form the lattice $Y$. The attractor
(resp. repellent) of $2\rho(\BC^*)$ to a fixed point $\mu$ is the orbit
$N'(\CK)^\varsigma\cdot\mu=:S_\mu$ (resp. $N'_-(\CK)^\varsigma\cdot\mu=:T_\mu$).
According to~\cite[3.3.2]{ti},
$\on{Gr}=\bigsqcup_{\mu\in Y}S_\mu=\bigsqcup_{\mu\in Y}T_\mu$.

\lem{mivi}
a) The closure $\overline{T}_\mu=\bigcup_{\nu\geq\mu}T_\nu$;

b) The closure $\overline{S}_\mu=\bigcup_{\nu\leq\mu}S_\nu$;

c) There is an isomorphism $\CF^\alpha_0\simeq S_0\cap\overline{T}_{-\alpha}$.
\elem

\prf
a) and b): same as~\cite[Proposition~3.1]{mv}.
c): same as~\cite[Theorem~2.7]{bfk}.
\epr

\lem{zhuzhu}
$\dim\CF^\alpha_0=|\alpha|$.
\elem

\prf
Same as~\cite[Theorem~3.2]{mv}, provided we know the dimensions of
$G'(\CO)^\varsigma$-orbits in the twisted Grassmannian: $\dim\on{Gr}^\eta=2|\eta|$
for $\eta\in Y^+$, according to e.g.~\cite[Corollary~2.10]{r}.
\epr

\cor{equidim}
Any fiber of $\pi:\ Z^\alpha\to\BA^\alpha$ is equidimensional of dimension
$|\alpha|$.
\ecor

\prf
Factorization.
\epr

\ssec{conn}{Irreducibility}
We consider the open subscheme $\overset{\circ}{Z}{}^\alpha:=
(\overset{\circ}{Z}{}^{a(\alpha)}_{\fg'})^\varsigma\subset Z^\alpha$ formed by the
based twisted maps (as opposed to quasimaps). The smoothness of
$\overset{\circ}{Z}{}^{a(\alpha)}_{\fg'}$ implies the smoothness of
$\overset{\circ}{Z}{}^\alpha$.

\prop{con}
$\overset{\circ}{Z}{}^\alpha$ is connected.
\eprop

\prf
We argue as in~\cite[Proposition~2.25]{bfg}. By induction in $\alpha$ and
factorization, if there are more than one connected components, we may
(and will) suppose that one of them, say $K'$, has the property
$\pi(K')\subset\Delta$ where $\Delta\subset\BA^\alpha$ is the main quasidiagonal.
By~\refc{equidim}, $\dim K'\leq|\alpha|+1$. By the same~\refc{equidim}, there
is another component $K$ such that $\pi(K)=\BA^\alpha$, and $\dim K=2|\alpha|$.
In the case $|\alpha|=1$ (i.e. $\alpha$ is a simple root of
$(\check G,\check T)$) we are reduced
to one of the two situations: a) $\fg'={\mathfrak{sl}}_2$, and the degree
$a(\alpha)$ is $d$ (long root $\alpha$);
b) $\fg'={\mathfrak{sl}}_2^{\oplus d}$, and the degree $a(\alpha)$
is 1 along each factor (short root $\alpha$). In both situations one checks
immediately $Z^\alpha\simeq\BA^2$.
So we may assume $|\alpha|>1$, and hence $\dim K>\dim K'$. This inequality
will lead to a contradiction. For $\phi\in K$ we
have $\dim K=\dim T_\phi\overset{\circ}{Z}{}^\alpha$. We have
$T_\phi\overset{\circ}{Z}{}^\alpha=H^0(\bC,\phi^*\CT\CB_{\fg'}(-\infty_\bC))^\Xi$
where $\CT\CB_{\fg'}$ stands for the tangent bundle of the flag variety
$\CB_{\fg'}=G'/B'$.
Since $\CT\CB_{\fg'}$ is generated by the global sections,
$H^0(\bC,\phi^*\CT\CB_{\fg'}(-\infty_\bC))=0$, and
$\dim T_\phi\overset{\circ}{Z}{}^\alpha$ can be computed as the invariant
part of the equivariant Euler characteristic of $\phi^*\CT\CB_{\fg'}(-\infty_\bC)$.
By the Atiyah-Singer equivariant index formula~\cite{as},
$\chi(\varsigma,\bC,\phi^*\CT\CB_{\fg'}(-\infty_\bC))$ is independent of $\phi$,
i.e. is the same for $\phi\in K$ and $\phi'\in K'$. Hence $\dim K=\dim K'$,
a contradiction.
\epr

\cor{irr}
$Z^\alpha$ is irreducible.
\ecor

\prf
We have to prove that $Z^\alpha$ is the closure of $\overset{\circ}{Z}{}^\alpha$.
The stratification $Z^{a(\alpha)}_{\fg'}=
\bigsqcup_{\Lambda_+\ni\gamma\leq a(\alpha)}\overset{\circ}{Z}{}^{\gamma}_{\fg'}\times
(\bC-\infty)^{\alpha-\gamma}$ induces the stratification
$Z^\alpha=\bigsqcup_{\beta\leq\alpha}\overset{\circ}{Z}{}^\beta\times\BA^{\alpha-\beta}$.
We argue as in~\cite[Theorem~10.2]{bfg}. It suffices to prove that
$(\phi,\unl{z})\in\overset{\circ}{Z}{}^\beta\times\BA^{\alpha-\beta}$ lies in the
closure of $\overset{\circ}{Z}{}^\alpha$ for $\unl{z}$ lying away from all
the quasidiagonals and distinct from $\pi(\phi)$. By factorization this
reduces to the case of simple $\alpha$. In this case $Z^\alpha\simeq\BA^2$
is irreducible, as was explained in the proof of~\refp{con}.
\epr

\ssec{normal}{Normality}
Recall that each $W$-orbit in $Y$ has a unique representative
$\eta$ such that
$a(\eta)\in X_*^+(T')$ is a dominant coweight. We call such $\eta$ dominant
as well, and we denote by $Y^+$ the cone of all dominant elements. Thus
$Y^+\iso Y/W\simeq G'(\CO)^\varsigma\backslash G'(\CK)^\varsigma/G'(\CO)^\varsigma$.
We define the congruence subgroup $\bK_{-1}\subset G'(\CK)^\varsigma$ as the
kernel of the evaluation morphism $ev:\ G'(\BC[\st^{-1}])^\varsigma\to(G')^\sigma$.
Given $\eta\in Y^+$ we consider the orbit
$\CW_\eta:=\bK_{-1}\cdot\eta\subset\on{Gr}$. For $\lambda\geq\eta\in Y^+$ we
define the {\em transversal slice} $\overline\CW{}^\lambda_\eta$ as the
intersection $\overline{\on{Gr}}{}^\lambda\cap\CW_\eta$. It follows
from~\cite[Theorem~8.4]{pr} that $\overline\CW{}^\lambda_\eta$ is normal with
rational singularities.

\prop{norm}
$Z^\alpha$ is normal.
\eprop

\prf
As in~\cite[Theorem~2.8]{bf11} we construct a $T\times\BG_m$-equivariant
morphism $s^\lambda_\eta:\ \overline\CW{}^\lambda_\eta\to Z^\alpha$ for
$\alpha=\lambda-\eta$ (note that in the present nonsimply laced situtation,
for the longest element $w_0\in W$ we have $-w_0\alpha=\alpha$). 
More precisely, the desired morphism is just the
restriction of the similar morphism of {\em loc. cit.} to $\varsigma$-fixed
points. Similarly to {\em loc. cit.} we show that $s^\lambda_\eta$ induces an
isomorphism $(s^\lambda_\eta)^*:\ \BC[Z^\alpha]\to\BC[\overline\CW{}^\lambda_\eta]$ on
functions of degree less than or equal to $n\in\BN$ (with respect to the
action of $\BG_m$), provided $\eta$ is big enough. Now one deduces the
normality of $Z^\alpha$ from normality of $\overline\CW{}^\lambda_\eta$ as
in~\cite[Corollary~2.10]{bf11}.
\epr

\ssec{bound}{The boundary of $Z^\alpha$}
Recall the stratification
$Z^\alpha=\bigsqcup_{\beta\leq\alpha}\overset{\circ}{Z}{}^\beta\times\BA^{\alpha-\beta}$.
The closure of the stratum $\overset{\circ}{Z}{}^{\alpha-\gamma}\times\BA^{\gamma}$
is denoted $\partial_\gamma Z^\alpha$. The union
$\bigcup_{i\in I}\partial_{\alpha_i}Z^\alpha$ is denoted $\partial_1Z^\alpha$ and
is called the boundary of $Z^\alpha$. More generally, the union
$\bigcup_{|\gamma|\geq n}\partial_\gamma Z^\alpha$ is denoted $\partial_nZ^\alpha$
(with the reduced closed subscheme structure). The open subscheme
$Z^\alpha\setminus\partial_2Z^\alpha$ is denoted $\overset{\bullet}{Z}{}^\alpha$.
By factorization and the calculations for $|\alpha|=1$ (proof of~\refp{con}),
$\overset{\bullet}{Z}{}^\alpha$ is smooth. We are going to prove that
$\partial_1Z^\alpha\subset Z^\alpha$ with the reduced closed subscheme structure
is a Cartier divisor.
Recall the function $F_{a(\alpha)}$ on $Z^{a(\alpha)}_{\fg'}$ constructed
in~\cite[Section~4]{bf11}.

\prop{sqrt}
a) There is a function $F_\alpha\in\BC[Z^\alpha]$ such that
$F^{d}_\alpha=F_{a(\alpha)}|_{Z^\alpha}$.

b) $F_\alpha$ is an equation of $\partial_1Z^\alpha\subset Z^\alpha$.
\eprop

\prf
Let us denote $F_{a(\alpha)}|_{Z^\alpha}$ by $f_\alpha$ for short.
Recall that $F_{a(\alpha)}$ has simple zeroes at any boundary component of
$Z^{a(\alpha)}_{\fg'}$~\cite[Lemma~4.2]{bf11}. We first prove
that $f_\alpha$ vanishes to the order exactly $d$ at any boundary
component $\partial_{\alpha_i}Z^\alpha,\ i\in I$. We start with $i\in I_0$
(notations of~\refss{twist z}, a long simple root of $(\check G,\check T)$,
i.e. a $\Xi$-fixed point, say $i'$, in $I'$). The corresponding simple
coroot of $(G',T')$ will be denoted by $\alpha'_{i'}$.
Since $Z^\alpha$ is smooth at the generic point of $\partial_{\alpha_i}Z^\alpha$,
and $Z^{a(\alpha)}_{\fg'}$ is smooth at the generic point of
$\partial_{\alpha'_{i'}}Z^{a(\alpha)}_{\fg'}$, and set-theoretically
$\partial_{\alpha_i}Z^\alpha=Z^\alpha\cap\partial_{\alpha'_{i'}}Z^{a(\alpha)}_{\fg'}$,
we have to check that the multiplicity of intersection of
$Z^\alpha$ with $\partial_{\alpha'_{i'}}Z^{a(\alpha)}_{\fg'}$ is generically equal to
$d$.
By factorization, we are reduced to the case $\fg'={\mathfrak{sl}}_2,\
a(\alpha)=d$. Then $Z^{a(\alpha)}_{\fg'}$ is the moduli space of
pairs of polynomials $(P(\bt),Q(\bt))$, $P$ monic of degree $d$,
$Q$ of degree less than $d$. Furthermore, $F_{a(\alpha)}$ is
the resultant $\on{Res}(P,Q)$. For the sake of definiteness, let
$d=3$. Then $Z^{a(\alpha)}_{\fg'}=\{(P=\bt^3+a_2\bt^2+a_1\bt+a_0,
Q=b_2\bt^2+b_1\bt+b_0)\}$, and $Z^\alpha$ is cut out by the equations
$a_2=a_1=b_2=b_1=0$. Then we have $\on{Res}(P,Q)|_{Z^\alpha}=b_0^3$. This takes
care of the case of a long simple root $\alpha_i$.

Now let $i\in I_1$ be a short simple root of $(\check G,\check T)$
corresponding to a free $\Xi$-orbit, say $i',i'',i'''$, in $I'$ (again,
for the sake of definiteness, we take $d=3$). Then $i',i'',i'''$
are all disjoint in the Dynkin diagram of $\fg'$, and the
intersection
$\partial_{\alpha'_{i'}}Z^{a(\alpha)}_{\fg'}\cap
\partial_{\alpha'_{i''}}Z^{a(\alpha)}_{\fg'}
\cap\partial_{\alpha'_{i'''}}Z^{a(\alpha)}_{\fg'}$
is generically transversal. Moreover, each of
$\partial_{\alpha'_{i'}}Z^{a(\alpha)}_{\fg'},
\partial_{\alpha'_{i''}}Z^{a(\alpha)}_{\fg'},
\partial_{\alpha'_{i'''}}Z^{a(\alpha)}_{\fg'}$
is generically transversal to
$Z^\alpha\subset Z^{a(\alpha)}_{\fg'}$, and generically
$\partial_{\alpha_i}Z^\alpha=Z^\alpha\cap
\partial_{\alpha'_{i'}}Z^{a(\alpha)}_{\fg'}=Z^\alpha\cap
\partial_{\alpha'_{i''}}Z^{a(\alpha)}_{\fg'}=
Z^\alpha\cap\partial_{\alpha'_{i'''}}Z^{a(\alpha)}_{\fg'}=
Z^\alpha\cap
\partial_{\alpha'_{i'}}Z^{a(\alpha)}_{\fg'}\cap
\partial_{\alpha'_{i''}}Z^{a(\alpha)}_{\fg'}
\cap\partial_{\alpha'_{i'''}}Z^{a(\alpha)}_{\fg'}$. This takes care of the case of
a short simple root $\alpha_i$.

We have $f_\alpha:\ \overset{\circ}{Z}{}^\alpha\to\BC^*$, and
$\sqrt[d]{f_\alpha}$ is well defined on an unramified Galois
covering $\widetilde{Z}\to\overset{\circ}{Z}{}^\alpha$ with Galois group $\Xi$.
To show the existence of $F_\alpha$ we have to prove that this covering splits,
i.e. the corresponding class in $H^1(\overset{\circ}{Z}{}^\alpha,\Xi)$ vanishes.
This is the subject of the following

\lem{tuesday}
There is a regular nonvanishing function
$F_\alpha\in\BC[\overset{\circ}{Z}{}^\alpha]$ such that
$F_\alpha^{d}=f_\alpha$.
\elem

\prf
Given a positive coroot $\alpha'$ of $G'$ we consider the moduli stack
$A_{\fg'}^{\alpha'}$ of $B'$-bundles over $C$ equipped with trivialization at
$\infty\in C$, such that the induced $T'$-bundle has degree $\alpha'$.
One can check that in case $\alpha'$ is dominant (as a coweight of $G'$)
$A_{\fg'}^{\alpha'}\simeq\BA^{2|\alpha'|}$. In general, $A_{\fg'}^{\alpha'}$
is a quotient of an affine space, and the automorphism groups of all
points are unipotent. The natural
morphism $\oZ^{\alpha'}_{\fg'}\to A_{\fg'}^{\alpha'}$ is an affine open
embedding with the image formed by all the $B'$-bundles $\phi_{B'}$ such that
the induced $G'$-bundle $\phi_{G'}$ is trivial. The complement divisor
$\fD_{\fg'}^{\alpha'}=A_{\fg'}^{\alpha'}\setminus\oZ^{\alpha'}_{\fg'}$ is
irreducible, and $F_{\alpha'}^{-1}$ extends to a regular function $F'_{\alpha'}$ on $A_{\fg'}^{\alpha'}$ vanishing to the order 1 along $\fD_{\fg'}^{\alpha'}$.

In case $\alpha'=a(\alpha)$, the automorphism
$\varsigma:\ \oZ^{a(\alpha)}_{\fg'}\to\oZ^{a(\alpha)}_{\fg'}$ extends to the
same named automorphism
$\varsigma:\ A^{a(\alpha)}_{\fg'}\to A^{a(\alpha)}_{\fg'}$, and we denote
the connected component of the fixed point stack
$(A^{a(\alpha)}_{\fg'})^\varsigma$ containing
$\oZ^\alpha$ by $A^\alpha$. One can check that in the appropriate coordinates
of the covering affine space $\BA^k\to A_{\fg'}^{a(\alpha)}$ the automorphism
$\varsigma$ is linear, so that $A^\alpha$ is also a quotient stack of an affine
space, and the automorphism groups of all points are unipotent as well. 
We denote the restriction
$F'_{a(\alpha)}|_{A^\alpha}$ by $f'_\alpha$ for short. The same argument as
in the first part of the proof of the proposition shows that $f'_\alpha$
vanishes to the order exactly $d$ at 
the complementary divisor $\fD^\alpha:=A^\alpha\setminus\oZ^\alpha$.

Finally, since $\pi_1(A^\alpha)=H_1(A^\alpha,\BZ)=H_1(A^\alpha,\Xi)=0$, the vanishing of the class in
$H^1(\oZ^\alpha,\Xi)$
associated to $\sqrt[d]{f'_\alpha}$ follows by excision from the above local
computations around $\fD^\alpha$.
\epr


So $F_\alpha$ is well defined on $\overset{\circ}{Z}{}^\alpha$, and extends
by zero through the generic points of the boundary divisor components
$\partial_{\alpha_i}Z^\alpha$. Hence it is defined off codimension 2, and extends
to the whole of $Z^\alpha$ by normality of $Z^\alpha$.

It remains to prove b), that is to check that the zero-subscheme of $F_\alpha$
is reduced. In other words, given $f\in\BC[Z^\alpha]$ vanishing at the
boundary $\partial_1Z^\alpha$ we have to check that $f$ is divisible by $F_\alpha$.
The rational function $f/F_\alpha$ is regular at the generic points of all
the boundary divisor components, so it is regular due to normality of
$Z^\alpha$.
\epr

\prop{degr}
$F_\alpha$ is an eigenfunction of $T\times\BG_m$ with the eigencharacter
$q^{(\alpha,\alpha)/2}z^{\alpha^*}$ (notations of~\refss{fold}).
\eprop

\prf Follows immediately from~\cite[Proposition~4.4]{bf11} along with
an observation that $d\cdot(\alpha,\alpha)=(a(\alpha),a(\alpha))$.
\epr

\rem{dets}
The invertible function $F_{a(\alpha)}|_{\overset{\circ}{Z}{}^{a(\alpha)}_{\fg'}}$
is constructed in~\cite[Section~4]{bf11} as the ratio of two sections of the
determinant line bundle lifted from $\on{Bun}_{G'}(\bC)$ (the generator
of its Picard group). The action of
$\Xi$ on $G'$ gives rise to a group scheme $\CG$ over $\bC/\!\!/\Xi$ as
in~\cite[Example~(3)]{he}. We have a natural morphism
$\on{Bun}_\CG\to\on{Bun}_{G'}(\bC)$, and the inverse image of the
determinant line bundle on $\on{Bun}_{G'}(\bC)$ is the determinant line bundle
on $\on{Bun}_\CG$ ({\em not} its $d$-th power), as
follows from~\cite[Theorem~3]{he} and~\cite[10.a.1,~(10.7)]{pr}.
\erem

\ssec{symp}{Canonical class of $Z^\alpha$}
Recall the smooth open subset $\overset{\bullet}{Z}{}^\alpha\subset Z^\alpha$ 
(notations of~\refss{bound}).

\lem{pole}
a) The canonical line bundle 
$\overset{\bullet}{\omega}_\alpha$ of $\overset{\bullet}{Z}{}^\alpha$ is trivial.

b) The weight of its generating section with respect to the loop rotations
is $\frac{(a(\alpha),a(\alpha))}{2}+|a(\alpha)|$.
\elem

\prf
a) For a simple coroot $\alpha_i$, we have $Z^{\alpha_i}\simeq\BA^2$; let us 
introduce coordinates $w_i$ along $\BA^{\alpha_i}$, and $z_i$ along the fibers
of $\pi_{\alpha_i}\colon Z^{\alpha_i}\to\BA^{\alpha_i}$. For arbitrary 
$\alpha=\sum a_i\alpha_i$, the factorization gives rise to the rational
coordinates $w_{i,r},z_{i,r},\ i\in I,\ 1\leq r\leq a_i$ on 
$\unl{Z}^\alpha:=Z^\alpha\times_{\BA^\alpha}\BA^{|\alpha|}$. Here the coordinates
$z_{i,r}$ are defined away from the quasidiagonal divisor in $\BA^{|\alpha|}$
(see~\refss{twist z}). Let us orient the Dynkin graph of $\fg$ so that
the arrows go from $I_0$ to $I_1$.
We define a rational section of the canonical line bundle
of $\overset{\bullet}{\unl{Z}}{}^\alpha$ by the following formula $s_\alpha:=
\prod_{i\in I}^{1\leq r<s\leq a_i}(w_{i,r}-w_{i,s})^2\cdot
\prod_{i\Rightarrow j}^{1\leq r\leq a_i,\ 1\leq s\leq a_j}(w_{i,r}-w_{j,s}^d)^{-1}\cdot
\prod_{i\to j}^{1\leq r\leq a_i,\ 1\leq s\leq a_j}(w_{i,r}-w_{j,s})^{-1}\cdot
\prod_{i\in I}^{1\leq r\leq a_i}dw_{i,r}dz_{i,r}$. Here $i\Rightarrow j$ means this
arrow belongs to the orientation chosen above and $i\in I_0,\ j\in I_1$;
and $i\to j$ means this arrow belongs to the orientation and $i,j\in I_0$
or $i,j\in I_1$.
This section extends regularly through the generic points of the irreducible
components of the quasidiagonal divisor due to 
examples in~\refss{exa},~\refss{anexa}
(for a divisor $w_{i,r}=w_{j,s},\ i\Rightarrow j$); due to~\cite[5.6]{bdf}
(for a divisor $w_{i,r}=w_{j,s},\ i\to j$); due to~\cite[5.5]{bdf}
(for a divisor $w_{i,r}=w_{i,s}$). Moreover, this section is 
$\fS_\alpha$-invariant, hence it descends to a rational section 
$\overline{s}_\alpha$ of $\overset{\bullet}{\omega}_\alpha$ that is regular
nonvanishing at the generic points of the irreducible components of the
quasidiagonal divisor (again due to~\cite[5.5]{bdf}). 
Hence $\overline{s}_\alpha$ trivializes $\overset{\bullet}{\omega}_\alpha$.

b) The weights of $w_{i,r},z_{i,r},dw_{i,r},dz_{i,r}$ with respect to the loop 
rotations are all equal to 1 if $i\in I_1$, and to $d$ if $i\in I_0$.
The explicit formula for $\overline{s}_\alpha$ implies the desired result.
\epr

\ssec{rs}{Rational singularities}
\prop{ra}
$Z^\alpha$ is a Gorenstein (hence, Cohen-Macaulay) scheme with canonical
(hence rational) singularities.
\eprop

\prf
We follow closely the proof of~\cite[Proposition~5.1]{bf11}, and use freely the
notations thereof. There we have
considered the Kontsevich resolution $\pi:\ M^{a(\alpha)}_{\fg'}\to Z^{a(\alpha)}_{\fg'}$, and computed its discrepancy divisor. Now we consider
the (smooth) fixed point stack $(M^{a(\alpha)}_{\fg'})^\Xi$
(see~\cite[especially~Proposition~3.7]{ro} for the basics on fixed point stacks with
respect to the finite groups' actions); more precisely, its
irreducible component $M^\alpha$ which is the closure of
$\overset{\circ}{Z}{}^\alpha\subset\overset{\circ}{Z}{}^{a(\alpha)}_{\fg'}
\subset M^{a(\alpha)}_{\fg'}$.\footnote{It is easy to see that
$(M^{a(\alpha)}_{\fg'})^\Xi$ is actually a special case of the moduli
space of twisted stable maps defined in~\cite{AV}.}
Note that there are other irreducible components
of $(M^{a(\alpha)}_{\fg'})^\Xi$, e.g. the loop rotation invariant stable maps
$(M^{a(\alpha)}_{\check\fg})^{\BG_m}$ (recall that $\CB=\CB_{\fg'}^\sigma$ is isomorphic
to $\sigma$-fixed points in the flag variety of $\check\fg'$ since $\fg'$ is
simply laced and hence isomorphic to $\check\fg'$. Hence $\CB$ is
isomorphic to the flag variety $\CB_{\check\fg}$ of $\check\fg$, and
$a(\alpha)\in H_2(\CB,\BZ)=H_2(\CB_{\fg'},\BZ)^\sigma=X_*(T')^\sigma$).
In notations of~\cite[proof of~Proposition~5.1]{bf11} the latter
component consists of stable maps such that $C=C_h\cup C_v$ where $\deg C_h=(1,0)$, and $\phi(C_h\cap C_v)=(0,\fb_-)$. This component is isomorphic to
the substack of based stable maps in $\overline{M}_{0,1}(\CB,a(\alpha))$,
and has dimension $2|a(\alpha)|-2$. Note also that the fixed point stack
$(M^{a(\alpha)}_{\fg'})^\Xi$ is {\em not} a closed substack of
$M^{a(\alpha)}_{\fg'}$: the natural morphism $(M^{a(\alpha)}_{\fg'})^\Xi\to M^{a(\alpha)}_{\fg'}$ has finite fibers over the points with nontrivial
automorphisms.

The complement $M^\alpha\setminus\overset{\circ}{Z}{}^\alpha$ is a union of
smoooth irreducible divisors $D_{\beta'}$ numbered by all $\beta'\in\Lambda_+$ (notations of~\refss{QQ}) such that
$\sum_{\xi\in\Xi}\xi(\beta')\leq a(\alpha)$. The generic point of $D_{\beta'}$
parametrizes the pairs $(C,\phi)$ such that $C=C_h\cup C_v$, the degree of
$\phi|_{C_h}$ equals $(1,a(\alpha)-\sum_{\xi\in\Xi}\xi(\beta'))$, and $C_v$
consists of irreducible components $C_v^\xi,\ \xi\in\Xi,\ \deg C_v^\xi=(0,\xi(\beta'))$ ($\Xi$-invariance implies in particular that the set
of points $\{C_v^\xi\cap C_h\}_{\xi\in\Xi}\subset C_h\simeq\BP^1$ is
$\Xi$-invariant). Among those divisors, $D_{\beta'}$ for {\em simple} $\beta'$
project generically one-to-one onto the boundary divisors of $Z^\alpha$.
The remaining divisors are exceptional.

The discrepancy of $\pi:\ M^\alpha\to Z^\alpha$ equals
$\sum_{\beta':\ \sum_{\xi\in\Xi}\xi(\beta')\leq a(\alpha)}m_{\beta'}D_{\beta'}$, and we have to 
show $m_{\beta'}\geq0$. As in~\cite[Proposition~5.1]{bf11}, by factorization
it suffices to consider the components $D_{\beta'}$ such that
$\sum_{\xi\in\Xi}\xi(\beta')=a(\alpha)$. The fixed point stack $D_{\beta'}^{\BG_m}$
with respect to the action of the loop rotations contains all the pairs
$(C,\phi)$ such that $C$ consists of $2+d$ irreducible components
$C_h,\ C_v^0,\ C_v^\xi,\ \xi\in\Xi,\ \deg C_h=(1,0),\
\deg C_v^\xi=(0,\xi(\beta')),\
\deg C_v^0=(0,0)$, with the following intersection pattern. The horizontal
component $C_h$ intersects $C_v^0$ at the point $0\in C_h\simeq\BP^1$.
The component $C_v^\xi$ intersects only $C_v^0$, and $\Xi$ acts on $C$
preserving $C_h,C_v^0$, and permuting the components $C_v^\xi,\ \xi\in\Xi$.
Note that the codimension of $D_{\beta'}^{\BG_m}$ in $D_{\beta'}$ is one.

We will prove $m_{\beta'}=|\beta'|+\frac{(\beta',\beta')}{2}-2$
(cf.~\cite[Lemma~5.2]{bf11}).
We will distinguish between the following two cases: a) {\em invariant case},
when $\beta'$ is $\Xi$-fixed; then the group of automorphisms of generic point
of $D_{\beta'}^{\BG_m}$ is equal to $\Xi$; b) {\em noninvariant case}, when
$\beta'\ne\xi\beta'$ for a nontrivial element $\xi\in\Xi$; then the group
of automorphisms of generic point of $D_{\beta'}^{\BG_m}$ is trivial.

We first consider the noninvariant case. Let $(C,\phi)\in D_{\beta'}$ be a
general point, and let $p_\xi:=C_v^\xi\cap C_h$. Then the fiber of the normal
bundle $\CN_{D_{\beta'}/M^\alpha}$ at the point $(C,\phi)$ equals
$(\bigoplus_{\xi\in\Xi}T_{p_\xi}C_v^\xi\otimes T_{p_\xi}C_h)^\Xi$. As $p_\xi\in C_h$
tends to $0\in C_h$, this tends to the fiber of $\CN_{D_{\beta'}/M^\alpha}$ at a
point $(\ 'C,\phi')$ of $D_{\beta'}^{\BG_m}$ equal to
$(\bigoplus_{\xi\in\Xi}T_{p_\xi}\ 'C_v^\xi\otimes T_0C_h)^\Xi$ where $p_\xi$ is
the intersection point of the components $'C_v^\xi$ and $'C_v^0$.
The group $\BG_m$ acts on this fiber
via the character $q^{-1}$ (cf.~\cite[proof of~Lemma~5.2]{bf11}).
On the other hand, the fiber of $\CN_{D_{\beta'}^{\BG_m}/D_{\beta'}}$ at the point
$(\ 'C,\phi')$ equals $T_0C_v^0\otimes T_0C_h$, and $\BG_m$ acts on this fiber
via the character $q^{-1}$ as well. Finally, $T_{(\ 'C,\phi')}D_{\beta'}^{\BG_m}$ is
nothing but $\Xi$-invariants in the similar tangent space described in
{\em loc. cit.} From this description it follows that $\BG_m$ acts trivially
on these invariants. All in all, $\BG_m$ acts on $\det T_{(\ 'C,\phi')}M^\alpha$
via the character $q^{-2}$, and on the fiber of the canonical bundle
$\omega_{M^\alpha}$ at $(\ 'C,\phi')$ via the character $q^2$.
Now the same argument as in {\em loc. cit.} yields
$m_{\beta'}=|\beta'|+\frac{(\beta',\beta')}{2}-2$.

In the invariant case, due to the presence of the automorphism group $\Xi$,
repeating the above argument, we obtain that $\BG_m$ acts on the fiber
of $\CN_{D_{\beta'}/M^\alpha}$ at $(\ 'C,\phi')$ via the character $q^{-d}$,
and on the fiber of $\omega_{M^\alpha}$ at $(\ 'C,\phi')$ via the character
$q^{2d}$. From this we deduce again
$m_{\beta'}=|\beta'|+\frac{(\beta',\beta')}{2}-2$.

Now we finish the proof of the proposition the same way as
in~\cite[proof of Proposition~5.1]{bf11}.
\epr

\ssec{vanish}{Cohomology vanishing}
Recall the notations of~\refss{QQ}. We
will consider the global quasimaps' spaces $\QM^\alpha$, and the
corresponding ind-scheme $\fQ$. We will generalize
the results of~\cite[Section~3]{bf12} on cohomology of the line bundles
$\CO_{\check\lambda}$ to the twisted case. We denote by
$\widetilde{H}{}^n(\fQ,\CO_{\check\lambda})$ the subspace of $\BG_m$-finite
vectors in $H^n(\fQ,\CO_{\check\lambda})$. Finally, given $\check\lambda\in X$,
we define a cofinal subsystem $Y_+^{\check\lambda}\subset Y_+$ formed by
$\alpha$ such that $\alpha^*+\check\lambda$ is dominant.

\prop{3.2}
(1) For $n>0$ and $\alpha\in Y_+^{\check\lambda}$ we have
$H^n(\QM^\alpha,\CO_{\check\lambda})=0$.

(2) For $n>0$ and $\check\lambda\in X$ we have
$\widetilde{H}{}^n(\fQ,\CO_{\check\lambda})=0$.

(3) For $\check\lambda\not\in X^+$ we have
$\widetilde{H}{}^0(\fQ,\CO_{\check\lambda})=0$.
\eprop

\prf
(3) is clear, and (2) follows from (1). We prove (1).

We will use the self evident notation $\partial_{\alpha_i}\QM^\alpha$ for the
boundary divisors of $\QM^\alpha$. We consider a divisor
$\Delta:=\sum_{i\in I}\partial_{\alpha_i}\QM^\alpha$.
We introduce the open subvariety
$\overset{\circ}\QM{}^\alpha\subset\QM^\alpha$ formed by all the twisted quasimaps
without defect at $\infty\in\bC$, and the evaluation morphism
$ev_\infty:\ \overset{\circ}\QM{}^\alpha\to\CB=(G'/B')^\sigma$. It is a fibration
with the fibers isomorphic to $Z^\alpha$. We have
$ev_\infty^*\omega_\CB=\CO_{-2\check\rho}$. It follows from~\refl{pole} that
$K_{\overset{\circ}\QM{}^\alpha}+\Delta-ev_\infty^*K_\CB=0$ (here $K$ stands for the
canonical class). According to~\refp{ra}, $Z^\alpha$ is Gorenstein with
rational singularities; but
$\QM^\alpha$ is locally in \'etale topology isomorphic to $Z^\alpha\times\CB$,
hence $\QM^\alpha$ is Gorenstein with rational singularities as well.
We conclude that the canonical bundle $\omega^\alpha:=
\omega_{\QM^\alpha}\simeq\CO_{\QM^\alpha}(-\Delta)\otimes\CO_{-2\check\rho}$.
We have the following analogue of~\cite[Lemma~4]{bf12}:

\lem{omegalpha}
$\omega^\alpha\simeq\CO_{-\alpha^*-2\check\rho}$.
\elem

\prf
As in the proof of~\cite[Lemma~4]{bf12} we see that there is $\check\mu\in X$
such that $\omega^\alpha\simeq\CO_{\check\mu}$. We have to check
$\check\mu=-\alpha^*-2\check\rho$. We will do this on an open subvariety
$\overset{\bullet}\QM{}^\alpha\subset\QM^\alpha$ with the complement of
codimension two. Namely, $\overset{\bullet}\QM{}^\alpha$ is formed by all
the twisted quasimaps of defect at most a simple coroot $\alpha_i,\ i\in I$
(or no defect at all). Note that $\Delta\cap\overset{\bullet}\QM{}^\alpha$
is a disjoint union of smooth divisors
$\partial_{\alpha_i}\overset{\bullet}\QM{}^\alpha$. Moreover,
$\overset{\bullet}\QM{}^\alpha$ itself is smooth, and the Kontsevich resolution
$K^\alpha\to\QM^\alpha$ (cf. proof of~\refp{ra}) is an isomorphism over
$\overset{\bullet}\QM{}^\alpha$. Let us fix a quasimap without defect
$\phi\in\QM^{\alpha-\alpha_i}$, choose a representative $\tilde\alpha_i$ of
$\alpha_i$, and consider a map
$p:\ \bC\to\partial_{\alpha_i}\overset{\bullet}\QM{}^\alpha$ sending $\bt\in\bC$ to
$\phi(\sum_{r=1}^{d}\sigma^r\tilde\alpha_i\cdot\zeta^{-r}\bt)$ (twisting
$\phi$ by a defect in $\bC^{a(\alpha_i)}$). Clearly, if $i\in I_1$ ($\alpha_i$
is a short root of $(\check G,\check T)$), then $p$ is a closed embedding; and if
$i\in I_0$ ($\alpha_i$ a long root of $(\check G,\check T)$), then $p$ factors through
$\bC\to\bC/\!\!/\Xi\hookrightarrow\partial_{\alpha_i}\overset{\bullet}\QM{}^\alpha$.
We will denote the categorical quotient $\bC/\!\!/\Xi$ (a projective line)
by $\overline\bC$, and its closed embedding into
$\partial_{\alpha_i}\overset{\bullet}\QM{}^\alpha$ by $\overline p$.
In both cases, the image of $\bC$ in
$\partial_{\alpha_i}\overset{\bullet}\QM{}^\alpha$ will be denoted by
$C^\phi_i$. It is easy to see that $\deg\CO_{\check\omega_j}|_{C^\phi_i}=
\delta_{ij}=\langle\alpha_i,\check\omega_j\rangle$. Hence it remains to check
that $\deg(\omega^\alpha|_{C^\phi_i})=-\langle\alpha_i,\alpha^*+2\check\rho\rangle$.
To this end recall that
$\omega^\alpha\simeq\CO_{\QM^\alpha}(-\Delta)\otimes\CO_{-2\check\rho}$, and
the Kontsevich resolution $K^\alpha\to\QM^\alpha$ is an isomorphism over
$\overset{\bullet}\QM{}^\alpha$. Thus we have to compute the degree of the
normal line bundle $\CN_{\partial_{\alpha_i}\overset{\bullet}K{}^\alpha/K^\alpha}|_{C^\phi_i}$
restricted to $C^\phi_i$, and prove
$\deg\CN_{\partial_{\alpha_i}\overset{\bullet}K{}^\alpha/K^\alpha}|_{C^\phi_i}=
\langle\alpha_i,\alpha^*\rangle$.

We follow the argument of~\cite[proof of~Proposition~4.4]{fk}, and consider
first the case $i\in I_1$. The universal stable map $(\CalC,\varphi)$ over
$C^\phi_i\subset K^\alpha$ looks as follows. For $\bt\in\bC\setminus\{0,\infty\}$
(recall that $C^\phi_i\simeq\bC$) the curve $\CalC_\bt$ has components
$C_h=\bC,C_v^r,\ 1\leq r\leq d$, and
$\varphi_\bt|_{C_h}=(\on{id},\phi)$, while
$\deg(\varphi_\bt|_{C_v^r})=(0,\sigma^r\tilde\alpha_i)$. The intersection point
$C_v^r\cap C_h$ is $\zeta^{-r}\bt$.
For $\bt=0$ (resp. $\infty$), the curve $\CalC_\bt$ has components
$C_h=\bC,C_v^0,C_v^r,\ 1\leq r\leq d$, and
$\varphi_\bt|_{C_h}=(\on{id},\phi)$, while $\deg(\varphi_\bt|_{C_v^0}=(0,0)$, and
$\deg(\varphi_\bt|_{C_v^r})=(0,\sigma^r\tilde\alpha_i)$. The intersection points
of the components all lie on $C_v^0$, and $C_v^0\cap C_h=0$ (resp. $\infty$).

The description of the normal bundle
$\deg\CN_{\partial_{\alpha_i}\overset{\bullet}K{}^\alpha/K^\alpha}$ given in the
proof of~\refp{ra} implies
$\deg\CN_{\partial_{\alpha_i}\overset{\bullet}K{}^\alpha/K^\alpha}|_{C^\phi_i}=
2+\langle\alpha_i,\alpha^*-\alpha_i^*\rangle=\langle\alpha_i,\alpha^*\rangle$.
The argument in the case $i\in I_0$ is similar.
\epr

Returning to the proof of the Proposition, it is finished the same way
as the one of~\cite[Theorem~3.2]{bf12}.
\epr

\sec{fer}{Fermionic formula and $q$-Whittaker functions}

\ssec{fermi}{Fermionic formula}
Recall the setup of~\refss{fold}. In particular, an isomorphism
$\alpha\mapsto\alpha^*$ from the root lattice of $(\check{G},\check{T})$ to
the root lattice of $(G,T)$ defined in the basis of simple roots as follows:
$\alpha_i^*:=\check{\alpha}_i$ (the corresponding simple coroot).
For an element $\alpha$ of the root lattice of $(\check{G},\check{T})$,
we denote by $z^{\alpha^*}$ the corresponding character of $T$. As usually,
$q$ stands for the identity character of $\BG_m$, and $q_i=q^{d_i}$.
For $\gamma=\sum_{i\in I}c_i\alpha_i$, we set
$(q)_\gamma:=\prod_{i\in I}\prod_{s=1}^{c_i}(1-q_i^s)$.

According to~\cite[Theorem~3.1]{ferm}, the recurrence relations
\eq{feigin}
\calJ_\alpha=\sum_{0\leq\beta\leq\alpha}
\frac{q^{(\beta,\beta)/2}z^{\beta^*}}{(q)_{\alpha-\beta}}\calJ_\beta
\end{equation}
uniquely define a collection of rational functions $\calJ_\alpha,\ \alpha\geq0$,
on
$T\times\BG_m$, provided $\calJ_0=1$. Moreover, these functions are
nothing but
the Shapovalov scalar products of the weight components of the Whittaker
vectors in the universal Verma module over the corresponding quantum group.

\th{main}
$\calJ_\alpha$ equals the character of $T\times\Gm$-module $\BC[Z^\alpha]$.
\eth

\prf
We have to prove that the collection of characters of $T\times\Gm$-modules
$\BC[Z^\alpha]$ satisfies the recursion relation~\refe{feigin}. Given the
geometric preparations undertaken in~\refs{twi qua}, the proof is the same as
the one of~\cite[Theorem~1.5]{bf11}.
\epr

We organize all $\calJ_\alpha$ into a generating function
$J_\fg^{\on{twisted}}(z,x,q)=\sum_{\alpha\in\Lambda_+}x^\alpha\calJ_\alpha$, the
equivariant twisted $K$-theoretic $J$-function of $\CB_{\fg'}$. The same way
as~\cite[Corollaries~1.6,1.8]{bf11} follow
from~\cite[Theorem~1.5]{bf11},~\reft{main} implies the following

\cor{dva}
The equivariant twisted $K$-theoretic $J$-function $J_\fg^{\on{twisted}}$ of
$\CB_{\fg'}$
is equal to the Whittaker matrix coefficient of the universal Verma module
of $U_q(\check\fg)$; it is an eigen-function of the quantum difference Toda
integrable system associated with $\check\fg$. \qed
\ecor

\ssec{weyl}{Twisted Weyl modules and $q$-Whittaker functions}

The notions of the local (resp. global) Weyl modules over the twisted current
algebra $(\fg'[\st])^\varsigma$ were introduced in~\cite{fok}
(resp.~\cite[Section~9]{cik}). Recall the notations of~\refss{twist z}.
Given a dominant $G$-weight
$\check\lambda=\sum_{i\in I}\langle\alpha_i,\check\lambda\rangle\check\omega_i$
we define
$\BA^{\check\lambda}:=\prod_{i\in I_1}(\bC-\infty)^{(\langle\alpha_i,\check\lambda\rangle)}\times
\prod_{i\in I_0}((\bC-\infty)/(\bt\mapsto\zeta^{-1}\bt))^{(\langle\alpha_i,\check\lambda\rangle)}$.
The character of $\BC[\BA^{{\check\lambda}}]$ with respect to the natural action of
$\BC^*$ is equal to ${\prod\limits_{i\in I}\prod\limits_{r=1}^{\langle\alpha_i,\check\lambda\rangle}(1-q_i^r)^{-1}}$. According to~\cite[Section~9]{cik} there exists
an action of $\BC[\BA^{{\check\lambda}}]$ on the global twisted Weyl
$(\fg'[\st])^\varsigma$-module $\CW^{\on{twisted}}({\check\lambda})$ such that

1) This action commutes with $(G'[\st])^\varsigma\rtimes \BC^*$;

2) $\CW^{\on{twisted}}({\check\lambda})$ is finitely generated and free
over $\BC[\BA^{{\check\lambda}}]$.

3) The fiber of $\CW^{\on{twisted}}({\check\lambda})$ at ${\check\lambda}\cdot 0\in \BA^{{\check\lambda}}$ is the
local twisted Weyl module $D^{\on{twisted}}({\check\lambda})$ of~\cite{fok}.

Now the local twisted Weyl modules $D^{\on{twisted}}({\check\lambda})$
coincide by~\cite{fok} with the level one Demazure modules over
$\fg'[\st]^\varsigma\rtimes\BC^*$. And the characters of level one Demazure
modules over dual untwisted affine Lie algebras were proved in~\cite{i} to
coincide with the $q$-Hermite polynomials
$\hat{\Psi}'_{{\check\lambda}}(q,z)$ (see~\refss{13}).

On the other hand, recall $q$-Whittaker functions $\Psi_{{\check\lambda}}(q,z)$
and $\hat{\Psi}_{{\check\lambda}}(q,z):=\Psi_{{\check\lambda}}(q,z)\cdot{\prod\limits_{i\in I}
\prod\limits_{r=1}^{\langle\alpha_i,\check\lambda\rangle}(1-q_i^r)}$
of~\cite[Theorem~1.2]{bf12}. Given the geometric preparations undertaken
in~\refs{twi qua}, the following theorem is proved the same way
as~\cite[Theorem~1.3]{bf12}:
\th{13}
The characters of $T\times\BC^*$-modules $\CW^{\on{twisted}}({\check\lambda})$ and
$D^{\on{twisted}}({\check\lambda})$ are
given by the corresponding $q$-Whittaker functions:
$\chi(\CW^{\on{twisted}}({\check\lambda}))=\Psi_{{\check\lambda}}(q,z);\
\chi(D^{\on{twisted}}({\check\lambda}))=\hat{\Psi}_{{\check\lambda}}(q,z)$.
\qed \eth

Also, the same argument as the one for~\cite[Theorem~1.5]{bf12} establishes
the following version of the Borel-Weil theorem for the dual global and local
twisted Weyl modules:
\th{15}
There is a natural isomorphism $\Gamma((G'[[\st]]/T'\cdot U'_-[[\st]])^\varsigma,
\CO({\check\lambda}))\simeq\Gamma(\fQ,\CO({\check\lambda}))
\simeq \CW^{\on{twisted}}({\check\lambda})^{\vee}$. Similarly,
$\Gamma((G'[[\st]]/B'_-[[\st]])^\varsigma,\CO({\check\lambda}))\simeq D^{\on{twisted}}({\check\lambda})^{\vee}$.
\eth

\sec{nono}{Nontwisted non simply laced case}

\ssec{no rs}{Quasimaps: rational singularities}
Recall that $\fg$ is a non simply laced simple Lie algebra, and $Z^\alpha_\fg$ is
the corresponding zastava space.

\prop{51}
$Z^\alpha_\fg$ has rational singularities.
\eprop

\prf
We are going to apply~\cite[Corollary~7.7]{fh}.
Recall~\cite[Definition~3.7]{fh} that an effective divisor $\Delta$ is called
a {\em boundary} on a variety $X$ if $K_X+\Delta$ is a $\BQ$-Cartier divisor.
We will take $X=Z^\alpha_\fg$, and $\Delta=\sum_{i\in I}\partial_{\alpha_i}Z^\alpha_\fg$
(the sum of boundary divisors $\partial_{\alpha_i}Z^\alpha_\fg$ with multiplicity one).
Recall the symplectic form $\Omega$ on $\overset\circ{Z}{}^\alpha_\fg$ constructed
in~\cite{fkmm}, and let $\Lambda^{|\alpha|}\Omega$ be the corresponding
regular nonvanishing section of $\omega_{\overset\circ{Z}{}^\alpha_\fg}$. According
to~\cite{fkmm}, $\Lambda^{|\alpha|}\Omega$ has a pole of the first order
at each boundary divisor component
$\partial_{\alpha_i}Z^\alpha_\fg\subset\overset\bullet{Z}{}^\alpha_\fg$. Here
$\overset\bullet{Z}{}^\alpha_\fg\subset Z^\alpha_\fg$ is an open smooth
subvariety with codimension
2 complement formed by all the quasimaps with defect of degree at most a
simple coroot. Recall a function $F_\alpha\in\BC[Z^\alpha_\fg]$~\cite[4.1]{bf11}.
According to~\cite[Lemma~4.2]{bf11}, $F_\alpha$ has a zero of order
$d_i=\frac{(\alpha_i,\alpha_i)}{2}$ at $\partial_{\alpha_i}Z^\alpha_\fg$.
Hence $F_\alpha\Lambda^{|\alpha|}\Omega$ is a regular section of
$\omega_{\overset\bullet{Z}{}^\alpha_\fg}$ nonvanishing at the boundary divisors
$\partial_{\alpha_i}Z^\alpha_\fg$ for a short coroot $\alpha_i$, and with a zero
of order $d_i-1$ for a long coroot $\alpha_i$.
We conclude that $\omega_{\overset\bullet{Z}{}^\alpha_\fg}\simeq
\CO_{\overset\bullet{Z}{}^\alpha_\fg}(\sum_{i\in I}(d_i-1)
\partial_{\alpha_i}Z^\alpha_\fg)$, and $K_{\overset\bullet{Z}{}^\alpha_\fg}+
\sum_{i\in I}\partial_{\alpha_i}Z^\alpha_\fg$ is the divisor of $F_\alpha$.
So indeed $\sum_{i\in I}\partial_{\alpha_i}Z^\alpha_\fg$ is a boundary on $Z^\alpha_\fg$
in the sense of~\cite[Definition~3.7]{fh}.

Recall~\cite[Proof of Proposition~5.1]{bf11} the Kontsevich resolution
$\pi:\ M^\alpha\to Z^\alpha_\fg$. According to~\cite[Definition~3.8]{fh}, the
{\em log relative canonical divisor}
$K^\Delta_{M^\alpha/Z^\alpha_\fg}:=K_{M^\alpha}+\Delta_M-\pi^*(K_{Z^\alpha_\fg}+\Delta)$ where
$\Delta_M$ is the proper transform of $\Delta$ on $M^\alpha$.
According to~\cite[Corollary~7.7]{fh}, if $K^\Delta_{M^\alpha/Z^\alpha_\fg}$ is a sum
of exceptional divisors of $M^\alpha$ with positive multiplicities, then
$Z^\alpha_\fg$ has rational singularities. So we have to compute the multiplicities
in $K^\Delta_{M^\alpha/Z^\alpha_\fg}$. We use freely the notations
of~\cite[Proof of Proposition~5.1]{bf11}. As in {\em loc. cit.}, by
factorization it suffices to compute the single multiplicity $m_\alpha$ of
$D_\alpha$. In case $\alpha=\alpha_i$ is simple, we have $m_{\alpha_i}=0$ by
the definition of $K^\Delta_{M^\alpha/Z^\alpha_\fg}$ since $D_{\alpha_i}$ is not exceptional
(note that this zero multiplicity
is {\em not} given by the formula of~\cite[Lemma~5.2]{bf11}).
In case $\alpha$ is {\em not} simple, the divisor $D_\alpha$ {\em is}
exceptional, and the argument in the proof of~\cite[Lemma~5.2]{bf11} goes
through word for word, giving the result
$m_\alpha=|\alpha|+\frac{(\alpha,\alpha)}{2}-2>0$.
This completes the proof of the proposition.
\epr

\ssec{no vanish}{Quasimaps: cohomology vanishing}
In this Section we follow the notations of~\cite{bf12}. In particular, we
will consider the global quasimaps' spaces $\QM^\alpha_\fg$, and the
corresponding ind-scheme $\fQ_\fg$. We will generalize
the results of~\cite[Section~3]{bf12} on cohomology of the line bundles
$\CO(\check\lambda)$ to the case of non simply laced $G$.

\prop{32}
(1) For $n>0$ and $\alpha\in\Lambda_+^{\check\lambda}$ we have
$H^n(\QM^\alpha_\fg,\CO(\check\lambda))=0$.

(2) For $n>0$ and $\check\lambda\in\Lambda^\vee$ we have
$\widetilde{H}{}^n(\fQ_\fg,\CO(\check\lambda))=0$.

(3) For $\check\lambda\not\in\Lambda^\vee_+$ we have
$\widetilde{H}{}^0(\fQ_\fg,\CO(\check\lambda))=0$.
\eprop

\prf
(3) is clear, and (2) follows from (1). We prove (1).

We will use the self evident notation $\partial_{\alpha_i}\QM^\alpha_\fg$ for the
boundary divisors of $\QM^\alpha_\fg$. We define the boundary
$\Delta_Q:=\sum_{i\in I}\partial_{\alpha_i}\QM^\alpha_\fg$. Recall the open subvariety
$\overset{\circ}\QM{}^\alpha_\fg\subset\QM^\alpha_\fg$ formed by all the quasimaps
without defect at $\infty\in\bC$, and the evaluation morphism
$ev_\infty:\ \overset{\circ}\QM{}^\alpha_\fg\to\CB_\fg$. It is a fibration with
the fibers isomorphic to $Z^\alpha_\fg$. We have
$ev_\infty^*\omega_{\CB_\fg}=\CO(-2\check\rho)$. The proof of~\refp{51} implies
$K_{\overset{\circ}\QM{}^\alpha_\fg}+\Delta_Q-ev_\infty^*K_{\CB_\fg}=0$.

Now we have $\CO(K_{\overset{\circ}\QM{}^\alpha_\fg}+\Delta_Q)=
\CO(-\alpha^*-2\check\rho)$. In effect, the proof of~\cite[Lemma~4]{bf12}
goes through word for word: first it suffices to check the equality on the
open subvariety $\overset{\bullet}\QM{}^\alpha_\fg\subset\QM^\alpha_\fg$ formed
by all the quasimaps with defect at most a simple root since the complement
$\QM^\alpha_\fg\setminus\overset{\bullet}\QM{}^\alpha_\fg$ has codimension two.
Second, it suffices to calculate the degree of the normal bundle
$\CN_{\partial_{\alpha_i}\QM^\alpha_\fg/\overset{\bullet}\QM{}^\alpha_\fg}$ restricted to the
curve $C_i^\phi$ defined in {\em loc. cit.} Third, the equality
$\deg\CN_{\partial_{\alpha_i}\QM^\alpha_\fg/\overset{\bullet}\QM{}^\alpha_\fg}|_{C_i^\phi}=
\langle\alpha_i,\alpha^*+2\check\rho\rangle$ is proved
in~\cite[Proposition~4.4]{fk}.

Finally, for $\alpha\in\Lambda^{\check\lambda}_+$ the line bundle
$\CL=\CO(\check\lambda)\otimes\CO(-K_{\QM^\alpha_\fg}-\Delta_Q)$ on
$\QM^\alpha_\fg$ is very ample. The vanishing of
$H^{>0}(\QM^\alpha_\fg,\CO(\check\lambda))=H^{>0}(\QM^\alpha_\fg,\CL\otimes
\CO(K_{\QM^\alpha_\fg}+\Delta_Q))$ follows from~\cite[Theorem~2.42]{fu}
which in turn is an immediate corollary of~\cite[Corollary~1.3]{kss}.
\epr

\bigskip

\footnotesize{ {\bf A.B.}: Department of Mathematics, University of Toronto
and Perimeter Institute of Theoretical Physics, Waterloo, Ontario, Canada,
N2L 2Y5;\\ 
Skolkovo Institute of Science and Technology;\\
{\tt braval@math.toronto.edu}}

\footnotesize{
{\bf M.F.}: National Research University
Higher School of Economics, Russian Federation,\\
Department of Mathematics, 6 Usacheva st, Moscow 119048;\\
Skolkovo Institute of Science and Technology;\\
{\tt fnklberg@gmail.com}}


\begin{thebibliography}{99}

\bibitem{AV}
D.~Abramovich and A.~Vistoli,
{\em Compactifying the space of stable maps},
J. Amer. Math. Soc. {\bf 15} (2002), 27--75.

\bibitem{as} M.~F.~Atiyah, I.~M.~Singer, {\em The index of elliptic
operators: III}, Ann. of Math. {\bf 87} (1968), 546--604.


\bibitem{bfg} A.~Braverman, M.~Finkelberg and D.~Gaitsgory, {\em Uhlenbeck
spaces via affine Lie algebras}, Progress in Math. {\bf 244} (2006), 17--135.

\bibitem{bfk} A.~Braverman, M.~Finkelberg and D.~Kazhdan, {\em Affine
Gindikin-Karpelevich formula via Uhlenbeck spaces}, Springer Proceedings in
Math. {\bf 9} (2012), 17--29.

\bibitem{bf11} A.~Braverman, M.~Finkelberg, {\em Semi-infinite Schubert
varieties and quantum $K$-theory of flag manifolds},
J. Amer. Math. Soc. {\bf 27} (2014), no. 4, 1147--1168.

\bibitem{bf12} A.~Braverman, M.~Finkelberg, {\em Weyl modules and q-Whittaker
functions}, Mathematische Annalen {\bf 359}, no. 1 (2014), 45--59.

\bibitem{bdf} A.~Braverman, G.~Dobrovolska, M.~Finkelberg, 
{\em Gaiotto-Witten superpotential and Whittaker $D$-modules on monopoles},
Adv.\ Math.\ {\bf 300} (2016), 451--472. 

\bibitem{c} I.~Cherednik, {\em Whittaker limits of difference spherical
functions}, Int. Math. Res. Notices, no. 20 (2009), 3793--3842.



\bibitem{cfk} V.~Chari, G.~Fourier, T.~Khandai, {\em A categorical approach to
Weyl modules}, Transform. Groups {\bf 15} (2010), no. 3, 517--549.


\bibitem{cik} V.~Chari, B.~Ion, D.~Kus, {\em Weyl modules for the hyperspecial
current algebra}, Int. Math. Res. Not. IMRN (2015), no. 15, 6470--6515. 

\bibitem{Et} P.~Etingof, {\em Whittaker functions on quantum groups and
$q$-deformed Toda operators}, Amer. Math. Soc. Transl. Ser. 2 {\bf
194} (1999), 9--25.


\bibitem{fh} T.~de~Fernex, C.~Hacon, {\em Singularities on normal varieties},
Compos. Math. {\bf 145} (2009), 393--414.

\bibitem{fm} B.~Feigin, M.~Finkelberg, I.~Mirkovi\'c, A.~Kuznetsov,
{\em Semi-infinite flags. II. Local and global Intersection Cohomology
of Quasimaps' spaces},  Differential topology, infinite-dimensional
Lie algebras, and applications, Amer. Math. Soc. Transl. Ser. 2, {\bf 194},
Amer. Math. Soc., Providence, RI (1999), 113--148.

\bibitem{ferm} B.~Feigin, E.~Feigin, M.~Jimbo, T.~Miwa, E.~Mukhin,
{\em Fermionic formulas for eigenfunctions of the difference Toda Hamiltonian},
Lett. Math. Phys. {\bf 88} (2009), 39--77.

\bibitem{fk} M.~Finkelberg, A.~Kuznetsov, {\em Global Intersection Cohomology
of Quasimaps' Spaces}, Int. Math. Res. Notices, no. 7 (1997), 301--328.

\bibitem{fkmm} M.~Finkelberg, A.~Kuznetsov, N.~Markarian, I.~Mirkovi\'c,
{\em A note on a symplectic structure on the Space of $G$-monopoles},
Commun. Math. Phys. {\bf 201} (1999), 411--421.
{\em Erratum}, Commun. Math. Phys. {\bf 334} (2015), 1153--1155;
arXiv:math/9803124, v6.

\bibitem{fok} G.~Fourier and  D.~Kus,
{\em Demazure and Weyl modules: The twisted current case},
Trans. Amer. Math. Soc. {\bf 365} (2013), no. 11, 6037--6064.

\bibitem{fu} O.~Fujino, {\em Introduction to the log minimal model program
for log canonical pairs}, draft available at
http://www.math.kyoto-u.ac.jp/~fujino/MMP21-s.pdf.



\bibitem{h} M.~Haiman, {\em Cherednik algebras, Macdonald polynomials and
combinatorics}, ICM 2006 Proceedings {\bf 3}, European Math. Soc. (2006),
843--872.

\bibitem{he} J.~Heinloth, {\em Uniformization of $\CG$-bundles}, Math. Ann.
{\bf 347} (2010), 499--528.

\bibitem{i}
B.~Ion, {\em Nonsymmetric Macdonald polynomials and Demazure characters},
Duke Math. J. {\bf 116} (2003), no. 2, 299--318.



\bibitem{kss} S.~Kov\'acs, K.~Schwede, K.~Smith, {\em The canonical sheaf
of Du Bois singularities}, Advances in Math. {\bf 224} (2010), 1618--1640.



\bibitem{mv} I.~Mirkovi\'c and K.~Vilonen, {\em Geometric Langlands duality
and representations of algebraic groups over commutative rings},
Annals of Math. (2) {\bf 166} (2007), 95--143.

\bibitem{pr} G.~Pappas, M.~Rapoport, {\em Twisted loop groups and their
affine flag varieties}, Advances in Math. {\bf 219} (2008), 118--198.

\bibitem{r} T.~Richarz, {\em Schubert varieties in twisted affine flag
varieties and local models}, J.~Algebra {\bf 375} (2013), 121--147.

\bibitem{ro} M.~Romagny, {\em Group actions on stacks and applications},
Michigan Math. J. {\bf 53} (2005), 209--236.

\bibitem{Sev} A.~Sevostyanov, {\em Regular Nilpotent Elements and Quantum
Groups}, Commun. Math. Phys. {\bf 204} (1999), 1--16.

\bibitem{ti} J.~Tits, {\em Reductive groups over local fields},
Automorphic forms, representations and L-functions (Proc. Sympos. Pure Math.,
Oregon State Univ., Corvallis, Ore., (1977), Part 1, pp. 29--69,
Proc. Sympos. Pure Math., XXXIII, Amer. Math. Soc., Providence, R.I., 1979.


\bibitem{z} X.~Zhu, {\em The geometric Satake correspondence for ramified
groups}, Ann. Sci. \'Ec. Norm. Sup\'er. (4) {\bf 48} (2015), no. 2, 409--451.

\end{thebibliography}
\end{document}